\documentclass[oneside,final]{amsart}
\usepackage{amssymb}
\usepackage{amsaddr}
\usepackage[utf8]{inputenc}
\usepackage[T1]{fontenc}
\usepackage[english]{babel}
\usepackage{csquotes}
\usepackage{comment}
\usepackage{xargs}
\usepackage{graphicx, color}
\setkeys{Gin}{draft=false}
\usepackage[figurewithin=none]{caption,subcaption} 
\graphicspath{{./}} \usepackage{enumitem}
\setlist[enumerate,1]{label={(\roman*)}}

\usepackage[notcite,notref,color,]{showkeys}

 \usepackage{hyperref}
 \hypersetup{
 	hidelinks,
   colorlinks  = true,    urlcolor    = blue,    linkcolor   = blue,    citecolor   = blue,      pdfauthor   = {Isabel Hubard, Elías Mochán, Antonio Montero },pdfsubject  = {Maniplexes},pdftitle    = {Symmetries of voltage operations},  }

\usepackage{tikz-cd}
\tikzcdset{row sep=tiny,arrows={shorten=-4pt}
}

\usepackage[loadshadowlibrary,shadow, colorinlistoftodos,textsize=tiny,obeyFinal]{todonotes}

\setuptodonotes{fancyline, backgroundcolor=gray!10,bordercolor=gray}

\newcommandx{\inline}[2][1=]{\todo[inline, #1]{#2}}
\newcounter{hw}
\newcommandx{\homework}[2][1=]{\todo[inline, caption={Homework \thehw} #1]{\stepcounter{hw} #2}}

 \newcommandx{\isaLong}[2][1=Long todo, usedefault]{\todo[inline, color=green!25, caption={\textbf{Isa:} #1}]{\textbf{Isa: #1}. #2}}
 \newcommandx{\eliasLong}[2][1=Long todo, usedefault]{\todo[inline, color=red!25,caption={\textbf{Elías:} #1}]{\textbf{Elías: #1}. #2}}
 \newcommandx{\teroLong}[2][1=Long todo, usedefault]{\todo[inline,color=blue!25, caption={\textbf{Tero:} #1}]{\textbf{Tero: #1}. #2} }

\makeatletter
    \providecommand\@dotsep{5}
\makeatother

\theoremstyle{plain}
\newtheorem{thm}{Theorem}[section]
\newtheorem{theorem}[thm]{Theorem}

\newtheorem{lemma}[thm]{Lemma}
\newtheorem{prop}[thm]{Proposition}
\newtheorem{proposition}[thm]{Proposition}
\newtheorem{coro}[thm]{Corollary}
\newtheorem{corollary}[thm]{Corollary}

\newtheorem*{conjecture*}{Conjecture}

\theoremstyle{definition}

\newtheorem{eg}[thm]{Example}

\theoremstyle{remark}

\newtheorem{remark}[thm]{Remark}
\newtheorem{question}[thm]{Question}

\numberwithin{equation}{section}

\usepackage{xargs}

\renewcommand{\leq}{\leqslant} \renewcommand{\geq}{\geqslant}
\renewcommand{\epsilon}{\varepsilon} \renewcommand{\subset}{\subseteq}  
\renewcommand{\{}{\lbrace}
\renewcommand{\}}{\rbrace}

\renewcommand{\bar}{\overline}

\newcommand{\bZ}{\mathbb{Z}}

\newcommand{\cC}{\mathcal{C}}
\newcommand{\cF}{\mathcal{F}}

\newcommand{\cK}{\mathcal{K}}
\newcommand{\cM}{\mathcal{M}}

\newcommand{\cP}{\mathcal{P}}

\newcommand{\cT}{\mathcal{T}}
\newcommand{\cU}{\mathcal{U}}
\newcommand{\cW}{\mathcal{W}}
\newcommand{\cX}{\mathcal{X}}
\newcommand{\cY}{\mathcal{Y}}
\newcommand{\cZ}{\mathcal{Z}}
\newcommand{\po}{\mathcal{P}}
\newcommand{\wW}{\mathcal{Z}}
\newcommand{\zZ}{\mathcal{W}}

\DeclareMathOperator{\aut}{\Gamma} 

\DeclareMathOperator{\stab}{Stab}

\DeclareMathOperator{\mon}{Mon}

\DeclareMathOperator{\fg}{\Pi}

\DeclareMathOperator{\oo}{\mathcal{O}}

\DeclareMathOperator{\med}{Med}
\DeclareMathOperator{\tr}{Tr}
\DeclareMathOperator{\norm}{Norm}
\DeclareMathOperator{\core}{Core}

\newcommand{\GXY}{\aut({\cX\ertimes\cY})}

\newcommandx{\Proot}[2][1=\cP, 2= \baseFlag, usedefault]{\left( #1,#2 \right)}

\newcommandx{\pth}[2][1=i, usedefault]{ {}^{#1} #2}

\newcommandx{\ertimes}[1][1=\eta, usedefault]{\rtimes_{#1}}
\newcommandx{\eltimes}[1][1=\eta, usedefault]{\rlimes_{#1}}

\newcommand{\1}{{\mathbf{1}}}
\newcommand{\2}{{\mathbf{2}}}

\newcommand{\XY}{\cX \ertimes \cY}

\newcommandx{\mlink}[6][1,3,5,6]{\arrow[dash, #1]{#2}[font=\tiny, description, #3]{#4}[outer sep=3pt, #5]{#6}}
\newcommandx{\mdart}[6][1,3,5,6]{\arrow[#1]{#2}[font=\tiny, description, #3]{#4}[outer sep=3pt, #5]{#6}}
\newcommandx{\msemi}[6][1,3,5,6]{\arrow[dash, #1]{#2}[font=\tiny, description, near start, #3]{#4}[near end, #5]{#6}}

\newcommand{\gen}[1]{\left\langle #1 \right\rangle}

\usepackage{dynkin-diagrams}
\newcommand\circleRoot[1]{ \draw[fill=white] (root #1) circle(3pt); \fill[black] (root #1) circle(1.5pt);}
\newcommandx{\coxeter}[5][3=, 4=, 5=, usedefault]{\begin{dynkinDiagram}[Coxeter, labels={#4}, #5]#1 #2 \foreach\r in {#3}{\circleRoot{\r}}
\end{dynkinDiagram}} 

\usepackage[sorting=nyt, backend=biber, style=numeric,defernumbers=true,url=false,eprint=false,giveninits=true,
maxbibnames=9,isbn=false,sortcites = true,citestyle=numeric-comp,
]{biblatex}
\begin{document}

\title{Symmetries of voltage operations on polytopes, maps and maniplexes}
\author{Isabel Hubard}
\address{Institute of Mathematics, National Autonomous University of Mexico (IM UNAM), 04510 Mexico City, Mexico}
\email{isahubard@im.unam.mx}
\author{Elías Mochán} 
\address{Department of Mathematics, Northeastern University, 02115 Boston, USA}
\email{j.mochanquesnel@northeastern.edu}
\author{Antonio Montero}
\address{Faculty of Mathematics and Physics, University of Ljubljana, 1000 Ljubljana, Slovenia}
\email{antonio.montero@fmf.uni-lj.si}

\begin{abstract}
Voltage operations extend traditional geometric and combinatorial operations (such as medial, truncation, prism, and pyramid over a polytope) to operations on maniplexes, maps, polytopes, and hypertopes.
In classical operations, the symmetries of the original object remain in the operated one, but sometimes additional symmetries are created; the same situation arises with voltage operations.
We characterise the automorphisms of the operated object that are derived from the original one and use this to bound the number of flag orbits (under the action of the automorphism group) of the operated object in terms of the original one.
The conditions under which the automorphism group of the original object is the same as the automorphism group of the operated object are given.
We also look at the cases where there is additional symmetry, which can be accurately described due to the symmetries of the operation itself.

\end{abstract}

\maketitle

\listoftodos\relax

\section{Introduction} \label{sec:intro}
Two of the more well-known operations on maps and polyhedra are the medial and the truncation. These are two of the Wythoffian operations (\cite[Section 5.7]{Coxeter_1973_RegularPolytopes}), which are defined on convex regular polyhedra as the convex hull of the orbit of point in the intersection of certain mirror reflections from the original polyhedra. 
These operations can be generalised to maps and abstract polyhedra,  even if they are not regular (see \cite{HubardOrbanicIvicWeiss_2009_MonodromyGroupsSelf, OrbanicPellicerWeiss_2010_MapOperations$k$}). 
They induce local transformations on the original polyhedron: the new faces always correspond to either the faces or the vertex-figures of the original.
In the medial, the new vertices are the old edges and the new edges are the old "corners" (vertex-face incident pairs). 
In the truncation, the new vertices are the old arcs and the new edges are the old edges combined with the old corners.

These operations have been well-understood in different contexts such as convex or skeletal polyhedra and maps on surfaces, both combinatorially and geometrically (see for example \cite{Coxeter_1973_RegularPolytopes, OrbanicPellicerWeiss_2010_MapOperations$k$}).
In particular, the symmetries of the medial or truncation of a map or polytope have been studied.
The automorphisms of a polyhedron $\po$ are automorphisms of its medial and truncated polyhedra. 
These two polyhedra may have additional automorphisms that are not derived from the automorphisms of $\po$ (see \cite{OrbanicPellicerWeiss_2010_MapOperations$k$}).
The extra symmetry of the medial of a polyhedron $\po$ is simple to explain: it is present if and only if $\po$ is self-dual. 
On the other hand, the extra symmetry of the truncation is more difficult to identify. 
Nevertheless, it is possible to find a family of polyhedra whose symmetries are the same as those of the original polyhedron. 
This family consists of polyhedra whose flags cannot be coloured with two colours in such a way that $2$-adjacent flags have different colours while $0$- and $1$- adjacent flags have the same colour. 
However, even if the flags of the polyhedra can be coloured in this way, it does not necessarily mean that the truncated has extra symmetry.

Although dissimilar, these two conditions have something in common: they enable us to identify polyhedra whose medial and/or truncation do not possess any additional symmetries.
We can deduce the number of orbits and their local arrangement in the medial and truncation of a polyhedron $\po$ if we are aware of the number of orbits and their local configuration under the group of automorphisms. 
More precisely: we can find $\po$ to be such that if $\cT$ is the {\em symmetry type graph} of $\po$ (the quotient of the flag graph of $\po$ by its automorphism group), then $\med(\cT)$ and $\tr(\cT)$ are the symmetry type graphs of $\med(\po)$ and $\tr(\po)$, respectively.

Voltage operations, first introduced in \cite{HubardMochanMontero_2023_VoltageOperationsManiplexes}, are a generalisation of classic geometric and topological operations on maps and polytopes (e.g. Wythoffian constructions, snub, dual, Petrie-dual, products of polytopes as described in \cite{GleasonHubard_2018_ProductsAbstractPolytopes}, and many more).
These operations may be a useful tool in addressing one of the fundamental questions of the theory of symmetric polytopes, as will be discussed later.

In a similar way to the case with the medial and the truncation, when we apply a voltage operation to a map or polytope, we usually obtain a new map or polytope with the same automorphism group, but with an increased number of flags. 
Furthermore, \cite[Theorem 5.1]{HubardMochanMontero_2023_VoltageOperationsManiplexes} states that the symmetry type graph of a map or polytope that has undergone a voltage operation can be expressed in terms of the symmetry type graph of the original object.

However, there is a caveat: just as with the medial and truncation, the resulting object may have more symmetries than the original one. 
The symmetry type graph of the resulting object with respect to the original group of automorphisms is revealed by \cite[Theorem 5.1]{HubardMochanMontero_2023_VoltageOperationsManiplexes}, however, it does not provide information about the symmetry type graph with respect to the full group (if it is larger).
Generally, the actual symmetry type graph 
is a quotient of the one we get from the theorem, and the answer is exact only when there are no additional symmetries.

The paper is organised as follows.
Section\nobreakspace \ref {sec:basics} gives a brief introduction to the concepts used throughout the paper.
Section\nobreakspace \ref {sec:voltageOperations} recalls the definitions and base properties of voltage operations given in \cite{HubardMochanMontero_2023_VoltageOperationsManiplexes}.
In Section\nobreakspace \ref {sec:autoX} we characterise the automorphisms of an operated object that come from automorphisms from the original object, while in Section\nobreakspace \ref {sec:Factions}
we explore when can we ensure that the automorphism group of the original object coincides with that of the operated one.
This exploration is done with the so-called $F$-actions, which are also introduced in Section\nobreakspace \ref {sec:Factions}.
Section\nobreakspace \ref {sec:autoY} investigates which symmetries of the operation always lift, and which ones lift for specific cases.
Finally, in Section\nobreakspace \ref {sec:concluding} we make some concluding remarks.

\section{Basic notions} \label{sec:basics}
We start by giving some of the definitions, basic notions and results that we use throughout the paper.
Many of these notions are treated in our previous paper \cite{HubardMochanMontero_2023_VoltageOperationsManiplexes}; we refer the reader to this manuscript for deeper detail.

This paper deals with voltage operations, which implies we shall deal with voltage graphs. 
Thus, we use a very general definition of graph: a \emph{graph} is a 4-tuple $X=(V(X),D(X),I,(\cdot)^{-1})$ where $V=V(X)$ is a non-empty set whose elements are called \emph{vertices}, $D=D(X)$ is a set whose elements are called \emph{darts}, $I:D\to V$ is a function, and $(\cdot)^{-1}:D\to D$ is an involution.
The vertex $I(d)$ is called the \emph{starting point of $d$}, while $T(D):=I(d^{-1})$ is the \emph{terminal point} of $d$. 
The dart $d^{-1}$ is the \emph{inverse of $d$}.
An \emph{edge} of $X$ is a set $\{d,d^{-1}\}$ and we denote the set of edges of $X$ by $E(X)$.
A \emph{semi-edge} is a dart $d$ that is its own inverse. 
If the dart $d$ is a semi-edge, we abuse notation and refer to the edge $\{d\}$ as a semi-edge as well.

A \emph{path} in a graph $X$ is a finite sequence of darts $W=d_1d_2\ldots d_k$ where the terminal point of $d_i$ is the starting point of $d_{i+1}$ for $i \in \{1,\dots,k-1\}$. 
The \emph{starting point of $W$} is $u=I(d_1)$, while the \emph{terminal point of $W$} is $v=T(d_k)$. In this case we we say that $W$ \emph{goes from $u$ to $v$} (symbolically $W:u\to v$).
If $W:u\to u$ for some vertex $u$ we say that it is \emph{closed} and \emph{based at $u$}. If a path is not closed we say that it is \emph{open}.

A \emph{premaniplex of rank $n$} or simply \emph{$n$-premaniplex} is a graph $\cX$ whose edges (and therefore, darts) are coloured with elements of the set $\{0,1,\ldots, n-1\}$ in such a way that every vertex is the starting point of exactly one dart of each colour, and if $i,j \in \left\{0, \dots, n-1 \right\} $ are such that $\left| i-j \right| \geq 2 $, then the alternating paths of length 4 with colours $i,j$ are closed.
We shall refer to the vertices of a premaniplex as \emph{flags}. 
We often abuse notation and use $\cX$ both for a premaniplex and for its set of flags, and use $x \in \cX$ to say that $x$ is a flag of $\cX$.
A connected premaniplex with no parallel edges and no semi-edges is called a {\em maniplex}.

Maniplexes generalise the flag-graph of maps (as 2-cellular embeddings of graphs) and geometric and abstract polytopes (in the sense of \cite{McMullenSchulte_2002_AbstractRegularPolytopes}).
Given a map or a polytope $\po$, one may define its flag-graph $\cX=\cX(\po)$ as the edge coloured graph whose vertices are the flags of $\po$, and the $i$-edges are between $i$-adjacent flags.
It is straightforward to see that these flag-graphs are maniplexes.
Moreover, the $i$-faces of $\po$ are in correspondence with the connected components $\cX_i$, resulting by removing all the $i$-edges of $\cX$. 
The incidence between the $i$-faces of $\po$ is given by the non-empty intersection of the corresponding connected components.
Even though given a maniplex it is not always the flag-graph of a polytope (see \cite{GarzaVargasHubard_2018_PolytopalityManiplexes}), whenever we start with a polytope we shall abuse notation and refer to the polytope and the maniplex (its flag-graph) as the same object.

Given a flag $x\in\cX$ and a colour $i\in\{0,1,\ldots,n-1\}$,  $\pth{x}$ denotes the unique dart of colour $i$
 that starts at $x$.
 The terminal point of $\pth{x}$ is denoted by $x^i$, and it is called the \emph{$i$-adjacent flag of $x$}.

The \emph{universal rank-$n$ Coxeter group} with string diagram is the group $\cC^{n}=\gen{r_{0}, \dots r_{n-1}}$ defined by the following relations:
\begin{equation}\label{eq:relsReg}
\begin{aligned}
r_{i}^{2} &= \epsilon && \text{for all } i \in \{0, \dots, n-1\}, \\
(r_{i}r_{j})^{2} &= \epsilon && \text{if } |i-j| \geq 2.
\end{aligned}
\end{equation}
When the rank is clear from the context, or irrelevant for our arguments, we omit the $n$ of $\cC^n$, and denote it only by $\cC$.
 Since most of our work in this paper is based on that of \cite{HubardMochanMontero_2023_VoltageOperationsManiplexes}, it is worth noticing that the group $\cC^n$ is denoted by $\mon(\cU)$ in \cite{HubardMochanMontero_2023_VoltageOperationsManiplexes}, where $\cU$ denotes the {\em universal maniplex of rank $n$}.

Equation\nobreakspace \textup {(\ref {eq:relsReg})} implies that, for any $n$-premaniplex $\cX$, the family of mappings $r_i: x \mapsto x^{i}$ ( $i \in \{0, \dots, n-1\}$)  defines an action of the group $\cC^{n}$ on the flags of $\cX$. 
This action is transitive if and only if $\cX$ is connected.
The \emph{monodromy group} (often also called the \emph{connection group}) of $\cX$ is the permutation group induced by this action.
Naturally, we denote this group by $\mon(\cX)$.

If $W:x\to x'$ is a path in $\cX$ whose darts have colours $i_1,i_2,\ldots,i_k$ (in that order), we  write $W$ as $\pth[w]{x}$, where $w=i_k\ldots i_2 i_1$. 
We may naturally associate the element $\omega=r_{i_k}\ldots r_{i_2} r_{i_1} \in \cC$ to the word $w$.
Let $w_1$ and $w_2$ be two words on $\{0,1,\ldots,n-1\}$ and let $\omega_1$ and $\omega_2$ be the respective elements of $\cC$ associated to those words in this way.
We say that the paths $\pth[w_1]{x_1}$ and $\pth[w_2]{x_2}$ are \emph{homotopic} if $x_1=x_2$ and $\omega_1=\omega_2$.
It is easy to see that homotopy is an equivalence relation.
We denote by $\pth[\omega]x$ the homotopy class of the path $\pth[w]{x}$. 
Note that any path in the homotopy class $\pth[\omega]{x}$ goes from $x$ to $\omega x$.
We will often not distinguish between a path and its homotopy class.

If a path $W_1$ ends at the starting point of another path $W_2$, these paths can be concatenated to get a new path $W_1 W_2$. 
These operation is well defined when we do it with homotopy classes: in fact $(\pth[\omega_1]{x})(\pth[\omega_2]{(\omega_1 x)}) = \pth[\omega_1 \omega_2]{x}$.
The set of homotopy classes of paths in a premaniplex $\cX$ with the concatenation as (partial) operation forms a groupoid $\fg(\cX)$, called \emph{the fundamental groupoid of $\cX$}. 
The subset of closed paths based at a flag $x$, denoted by $\fg^x(\cX)$, forms a group and it is called \emph{the fundamental group of $\cX$ based at $x$}.

Given two $n$-premaniplexes $\cX_{1}$ and $\cX_{2}$, a \emph{(premaniplex) homomorphism} from $\cX_{1}$ to $\cX_{2}$ is a function that preserves $i$-adjacencies, for $i\in\{0, \dots, n-1\}$.
A surjective homomorphism is called a \emph{covering}. 
If there is a covering from $\cX_{1}$ to $\cX_{2}$ we say that \emph{$\cX_{1}$ covers $\cX_{2}$}, and that $\cX_{2}$ is a \emph{quotient} of $\cX_{1}$. 
It is easy to see that if $\cX_{2}$ is connected, any homomorphism from $\cX_{1}$ to $\cX_{2}$ is a covering.

A bijective homomorphism is called an \emph{isomorphism}, and an isomorphism from $\cX$ to itself is called an \emph{automorphism}.
The set of all the automorphisms of $\cX$ forms the \emph{automorphism group} of $\cX$ and it is denoted by $\aut(\cX)$.
It is straightforward to see that automorphisms and monodromies commute. 
Moreover, every permutation of the flags of $\cX$ that commutes with every monodromy is an automorphism.

If $\aut(\cX)$ has exactly $k$ orbits on the set of flags of $\cX$, we say that $\cX$ is a {\em $k$-orbit} premaniplex.

Given a premaniplex $\cX$ and a group $\Gamma \leq \aut(\cX)$, we define $\cX/\Gamma$ in the following way: the flags are the orbits $\{x\Gamma  :  x \in V(\cX)\}$ and for $i \in \{0, \dots, n-1\}$ $(x\Gamma)^{i} = (x^{i})\Gamma$.
Since $(x^{i}) \gamma = (x \gamma)^{i}$ for every $\gamma \in \aut(\cX)$, these adjacencies are well-defined. 
It is straightforward that $\cX/\Gamma$ is an $n$-premaniplex and that the mapping $x \mapsto x\Gamma$ is a covering from $\cX$ to $\cX/\Gamma$.
The quotient $\cX/\Gamma$ is the \emph{symmetry type graph of $\cX$} with respect to $\Gamma$. 
Further, we refer to $\cX / \aut(\cX)$ simply as the \emph{symmetry type graph of $\cX$}.

\section{Voltage operations} \label{sec:voltageOperations}

In \cite{HubardMochanMontero_2023_VoltageOperationsManiplexes} we define a voltage operation for maniplexes (maps and polytopes), and a voltage operator. 
This section revisits such definitions as well as some useful results of these concepts.

A \emph{voltage assignment} on a premaniplex $\cY$ is  a function $\xi:\Pi(\cY)\to G$ where $G$ is some group, with the property that $\xi(W_{1}W_{2})=\xi(W_{2})\xi(W_{1})$ for any two paths $W_{1} ,W_{2}$ such that the product $W_{1}W_{2}$ exists. 
The element $\xi(W)\in G$ is called \emph{the voltage of $W$}.
Note that if $W=d_1 d_1\ldots d_k$ then $\xi(W) = \xi(d_k) \cdots \xi(d_2)\xi(d_1)$, so the voltage of a path is determined by the voltages of its darts.
Therefore, we often just define voltages by assigning an element $\xi(d)\in G$ to each dart $d$ in such a way that $\xi(d^{-1})=(\xi(d))^{-1}$.
Graphically, this is illustrated by writing the voltage next to an edge with an orientation that indicates which dart has the written voltage. If the voltage is an involution, no orientation is needed.
Voltage assignments can be defined for graphs in general (see~\cite{Gross_1974_VoltageGraphs}).

Let $\cX$ be an $n$-premaniplex and let $\cY$ be an $m$-premaniplex.
Consider a voltage assignment $\eta: \fg(\cY) \to \cC^n$.
We define the $m$-premaniplex $\cX \ertimes \cY$ in the following way: the flags of $\cX \ertimes \cY$ are those in the set $\cX \times \cY$ and, for each $i\in \{0,1,\dots, m\}$, there is an edge of colour $i$ from  $(x,y)$ to $\left(\eta(\pth{y})x,r_i y\right)$; that is
\[(x,y)^{i} =\left(\eta(\pth{y})x,r_i y\right). \]

We call the pair $(\cY,\eta)$ a \emph{voltage operator} (or an \emph{(n,m)-voltage operator} if we want to emphasise the ranks), and the mapping $\cX\mapsto \XY$, that turns $n$-premaniplexes into $m$-premaniplexes is the corresponding \emph{voltage operation} (or \emph{$(n,m)$-voltage operation}).

There are plenty of known operations that can be described as voltage operations. 
These include Wythoffian constructions, the prism and pyramid over a polytope or maniplex, the snub operation on a map, the orientable double cover (of a non-orientable polytope or maniplex), and many more (see \cite[Section 4]{HubardMochanMontero_2023_VoltageOperationsManiplexes}).

An important property of voltage operations is that they preserve the symmetries of the premaniplex we apply them to. 
Furthermore, if $\cX$ is a premaniplex, $\Gamma \leq \aut(\cX)$ and $(\cY,\eta)$ a voltage operator, then $(\XY)/\Gamma \cong (\cX/\Gamma)\ertimes \cY$ (see~\cite[Theorem 5.1]{HubardMochanMontero_2023_VoltageOperationsManiplexes}).
In fact this property characterises voltage operations (\cite[Theorem 7.3]{HubardMochanMontero_2023_VoltageOperationsManiplexes}).

However one must be careful: some operations, like the snub, seem to only preserve some of the symmetries of the operated object (e.g. the rotations). 
The reason is that when seen as voltage operations, some operations might give as a result a premaniplex with many connected components, and the ``missing'' symmetries appear as isomorphisms between those components.
For example, when seeing the snub as a voltage operator $(\cY,\eta)$ (as in \cite[Section 4]{HubardMochanMontero_2023_VoltageOperationsManiplexes})  for every orientable map $\cM$, the map $\cM\ertimes \cY$ consists of a left snub of $\cM$ and a disjoint right snub of $\cM$; and if $\cM$ admits an automorphism sending a flag to an adjacent one (e.g. if $\cM$ is regular), it appears as an isomorphism between these two copies.

While some results work for all voltage operations, others work only for those operations that {\em preserve connectivity}, that is, those where $\XY$ is connected whenever $\cX$ is connected. 
The following result characterises such operations.

\begin{prop}[{\cite[Corollary 3.9]{HubardMochanMontero_2023_VoltageOperationsManiplexes}}]\label{coro:connected}
    An $(n,m)$-voltage operator $(\cY,\eta)$ preserves connectivity if and only if $\cY$ is connected and $\eta(\fg^{y_0}(\cY)) = \cC^n$ for any $y_0\in \cY$.
\end{prop}

It is possible that different voltage assignments on $\cY$ give rise to the same voltage operation.
In this case we say that the corresponding voltage operators are \emph{equivalent} (the formal definition has some additional nuances, see~\cite[Proposition 7.2]{HubardMochanMontero_2023_VoltageOperationsManiplexes} for details).
Using known results from voltage assignments on graphs, one can show that any voltage operator is equivalent to one where $\cY$ has a spanning tree $T$ where all its darts have trivial voltage (i.e. $\xi(d)=1$ for every dart $d$ in $T$). 
\begin{remark}\label{remark:trivialspanningtree}
    Unless otherwise stated, we shall assume that $(\cY,\eta)$ has a spanning tree with trivial voltage.
\end{remark}

\section{Characterizing automorphisms induced by $\cX$.} \label{sec:autoX}

In \cite{HubardMochanMontero_2023_VoltageOperationsManiplexes}, we showed that each automorphism of $\cX$ can be regarded as an automorphism of $\cX \ertimes \cY$.
More precisely, we showed that every automorphism $\gamma \in \aut(\cX)$, induces an automorphism $\bar{\gamma} \in \aut(\cX \ertimes \cY)$ by $(x,y)\bar{\gamma} = (x\gamma, y)$. 
That is, the automorphisms of $\cX$ act on the $\cX$-coordinate of $\cX \ertimes \cY$, and fix the $\cY$-coordinate. 
On  the other hand, in~\cite[Theorem 5.1]{HubardMochanMontero_2023_VoltageOperationsManiplexes} we showed that the quotient $(\XY)/\aut(\cX)$ is isomorphic to $(\cX/\aut(\cX))\ertimes \cY$; that is, that the symmetry type graph of $\XY$ with respect to $\aut(\cX)$ is  $\cT\ertimes \cY$, where $\cT$ is the symmetry type graph of $\cX$.
If it happens that $\GXY=\aut(\cX)$, this implies  that the symmetry type graph of $\XY$ is precisely $\cT\ertimes \cY$.  Thus, under this hypothesis, if $\cX$ is a $k$-orbit premaniplex, then $\XY$ is $k|\cY|$-orbit one.
Unfortunately requiring that $\GXY=\aut(\cX)$ seems to be a rather strong hypothesis.

In this section we revisit the automorphisms of $\XY$ that come from symmetries of $\cX$, we also give a characterisation of them and use these results to bound the number of flag-orbits of $\XY$ whenever $(\cY,\eta)$ is a connectivity-preserving voltage operator.
Of course, not all automorphisms of $\XY$ need to be of this type. In Section\nobreakspace \ref {sec:ExtraSym} and Section\nobreakspace \ref {sec:autoY} we shall study automorphisms of $\XY$ that are not induced by automorphisms of $\cX$ (in the sense discussed above).

\begin{lemma}\label{lemma:autX}
    Let $\cX$ be a connected premaniplex and let $(\cY,\eta)$ be a voltage operator that preserves connectivity. 
    If $\gamma \in \aut(\cX \ertimes \cY)$  fixes the $\cY$-coordinate of every $(x,y) \in \cX \ertimes \cY$, then $\gamma$ can be regarded as an automorphism of $\cX$.
\end{lemma}

\begin{proof}
Since $\gamma$ fixes the $\cY$-coordinate, while (maybe) moving the $\cX$-coordinate, and $\cX$ is connected, we can write the $\cX$-coordinate of the image of $(x,y)$ under $\gamma$ as a monodromy that depends on both $x$ and $y$ acting on $x$.
Thus, we may define $\theta: \cX \ertimes \cY \to \cC$ such that $(x,y)\gamma = (\theta(x,y)x, y)$.
To show that $\gamma \in \aut(\cX)$ we first show that for $y,y' \in \cY$, $\theta(x,y)x=\theta(x,y')x$. 
Given $y, y' \in \cY$, let $\omega\in \cC$ be such that $y'=\omega y$ and $\eta(\pth[\omega]{y})=1$ (it exists since,  by~Proposition\nobreakspace \ref {coro:connected}, $\cY$ is connected, and by Remak\nobreakspace \ref {remark:trivialspanningtree} there is a spanning tree of trivial voltage). Then,
\begin{eqnarray*}
(\theta(x,y')x, y') &=& (x,y')\gamma = (x,\omega y)\gamma = (\eta(\pth[\omega]{y})x, wy)\gamma = (\omega (x,y)) \gamma = \\
&=& \omega ((x,y)\gamma) = \omega (\theta(x,y)x, y)  = (\eta(\pth[\omega]{y})\theta(x,y)x, \omega y) =  \\
&=& (\theta(x,y)x, y'). 
\end{eqnarray*}
Hence the action of $\theta$ on $x$ does not depend on the $\cY$-coordinate of $\cX \ertimes \cY$, and we may denote it as $\theta(x)$.

Let $\omega$ be a monodromy. Using the fact that $\gamma$ is an automorphism of $\cX\ertimes \cY$ we get that, on one hand:

\[
    (\omega(x,y))\gamma = (\eta(\pth[\omega]{y})x,\omega y)\gamma = (\theta(\eta(\pth[\omega]{y})x)\eta(\pth[\omega]{y})x, \omega y),
\]
and on the other:
\[
    (\omega(x,y))\gamma = \omega((x,y)\gamma) = \omega (\theta(x) x, y) = (\eta(\pth[\omega]{y})\theta(x) x, \omega y).
\]
Therefore:
\begin{equation}\label{eq:prop_theta}
    \theta(\eta(\pth[\omega]{y})x)\eta(\pth[\omega]{y})x = \eta(\pth[\omega]{y})\theta(x) x.
\end{equation}

Now, abusing notation, we define the action of $\gamma$ on $\cX$ as $x\gamma:= \theta(x)x$; we need to show that $\gamma$ acts as an automorphism of $\cX$.

Since $(\cY,\eta)$ preserves connectivity, for each $i\in\{0,\dots,n-1\}$ there is some closed path $\pth[\omega_i]y\in \fg^y(\cY)$ such that $\eta(\pth[\omega_i]{y}) = r_i$ (Proposition\nobreakspace \ref {coro:connected}). 
Therefore,
\begin{align*}
    (r_i x)\gamma &= (\eta(\pth[\omega_i]{y}) x) \gamma\\
    &= \theta(\eta(\pth[\omega_i]{y}) x) \eta(\pth[\omega_i]{y}) x \\
    &= \eta(\pth[\omega_i]{y})\theta(x) x & \textrm{ (because of Equation\nobreakspace \textup {(\ref {eq:prop_theta})})}\\
    &= r_i\theta(x) x\\
    &= r_i(x\gamma).
\end{align*}
This shows that $\gamma$ acts on $\cX$ as an automorphism.
\end{proof}

Lemma\nobreakspace \ref {lemma:autX} says that those automorphisms of $\XY$ that always preserve the $\cY$-coordinate are automorphisms of $\cX$. The following result shows that it is enough to see that the automorphism preserves the $\cY$-coordinate of a single flag.

\begin{prop}\label{lemma:autX-1}
    Let $\cX$ be a connected premaniplex and $(\cY,\eta)$ be a voltage operator that preserves connectivity. Let $\gamma\in\GXY$. If there exists
     $(x_0,y_0)\in\XY$ such that $(x_0,y_0)\gamma = (x_1,y_0)$ for some $x_1\in\cX$, then $\gamma$ may be regarded as an automorphism of $\cX$.
\end{prop}

\begin{proof}
    Let $(x_0,y_0)\gamma=(x_1,y_0)$ and let $(x,y)$ be any flag of $\XY$.
    We want to show that $\gamma$ fixes the $\cY$-coordinate of $(x,y)$ and use~Lemma\nobreakspace \ref {lemma:autX} to conclude the proof.
    
    Since $\XY$ is connected by hypothesis, there is some $\omega\in\cC^m$ such that $(x,y)=\omega(x_0,y_0)$ and in particular $y = \omega y_0$.
    Then,
    \begin{align*}
        (x,y)\gamma &= (\omega(x_0,y_0))\gamma\\
        &= \omega ((x_0,y_0)\gamma)\\
        &= \omega(x_1,y_0)\\
        &= (\eta(\pth[\omega]{y_0}) x_1, \omega y_0)\\
        &= (\eta(\pth[\omega]{y_0}) x_1, y).
    \end{align*}

This proves that $\gamma$ preserves the $\cY$-coordinate of every flag of $\XY$. Therefore, Lemma\nobreakspace \ref {lemma:autX} tells us that $\gamma$ may be regarded as an automorphism of $\cX$.
\end{proof}

Whenever the graph $\cY$ is finite, the above result helps us to bound the number of flag orbits (under the action of the full automorphism group) of $\cX \ertimes \cY$. 
The bound is in terms of the number of flag orbits of $\cX$, the cardinality of $\cY$ and the index $[\aut(\XY) : \aut(\cX)]$.
More precisely:

\begin{prop}\label{prop:cota_orbitas}
    Let $\cX$ be a $k$-orbit premaniplex and let $(\cY,\eta)$ be a voltage operator where $\cY$ is a finite premaniplex. 
    Let $\aut(\cX) \leq \Gamma \leq \GXY$.
The following hold: 
    \begin{enumerate}
        \item \label{item:orbits} 
        the number of flag-orbits of $\Gamma$ on $\XY$ is \[\frac{k |\cY|}{[\Gamma: \aut(\cX)]};\]
        \item \label{item:index} if $(\cY,\eta)$  preserves connectivity, then $1 \leq [\Gamma: \aut(\cX)] \leq |\cY|$;
        \item \label{item:indexXY} in particular, if $(\cY,\eta)$ preserves connectivity, and $[\Gamma: \aut(\cX)] = \ell$ then $\cX \ertimes \cY$ is a $\frac{k |\cY|}{t \ell}$-orbit premaniplex, where $t$ is an integer satisfying $1 \leq t \leq \frac{|\cY|}{\ell} $.
    \end{enumerate}
\end{prop}

\begin{proof}
For simplicity let $\Delta$ denote the group $\aut(\cX)$.
Consider the set $O=\{(x,y)\Delta : (x,y) \in \XY\}$ of $\Delta$-orbits of $\XY$. 
Notice that \[\{(x,y)\Delta : (x,y) \in \XY\} = \{(x\Delta,y) : (x,y) \in \XY\},\] which implies that $|O| = k|\cY|$.

Now observe that every $\Delta$-orbit of $\XY$ is contained in exactly one $\Gamma$-orbit of $\XY$. 
Moreover, if $T$ is a set of left coset representatives of $\Delta$ in $\Gamma$ then each $\Gamma$-orbit $(x,y)\Gamma$ can be partitioned as \[(x,y)\Gamma = \bigcup_{\alpha \in T} (x,y) \alpha \Delta.\] 
That is, each $\Gamma$-orbit of $\XY$ contains exactly $t = [\Gamma : \Delta]$ $\Delta$-orbits.
Therefore the number of $\Gamma$-orbits of $\XY$ is $\frac{|O|}{t} = \frac{k|\cY|}{t}$. 
This proves Item\nobreakspace \ref {item:orbits}.

To prove Item\nobreakspace \ref {item:index} just 
recall that 
if $(x_{1}, y_{1})$ and $(x_{2}, y_{2}) $ are both in the same $\Gamma$-orbit and satisfy that $y_{1} = y_{2}$, then they are actually in the same $\Delta$-orbit (see Proposition\nobreakspace \ref {lemma:autX-1}).
It follows that if two $\Delta$-orbits $(x_1,y_2) \Delta$, $(x_2, y_2) \Delta$ are contained in the same $\Gamma$-orbit and satisfy that $y_1=y_2$, then $(x_1,y_1) \Delta = (x_2, y_2) \Delta$.
Therefore, the number of $\Delta$-orbits contained in the same $\Gamma$-orbit is bounded by $|\cY|$ but  
we just saw that this number is $[\Gamma : \Delta]$.
Thus $1\leq  [\Gamma:\Delta] \leq |\cY|$.

Of course Item\nobreakspace \ref {item:orbits} holds for the group $\GXY$, which implies that the number of orbits of $\cX \ertimes \cY$ is 
\[\frac{k|\cY|}{[\GXY: \Delta ]} = \frac{k|\cY|}{[\GXY:\Gamma][\Gamma: \Delta ]} = \frac{k|\cY|}{[\GXY:\Gamma] \ell }.\]
Now, because of Item\nobreakspace \ref {item:index}, the number $t=[\GXY:\Gamma]$ satisfies that \[t \ell = [\GXY:\Gamma][\Gamma: \Delta ] = [\GXY:\Delta] \leq |\cY|, \]
which completes the proof of Item\nobreakspace \ref {item:indexXY}.
\end{proof}

Proposition\nobreakspace \ref {prop:cota_orbitas} generalises~\cite[Proposition 4.3]{OrbanicPellicerWeiss_2010_MapOperations$k$}, that says that the number of orbits of the truncation of a $k$-orbit map $\cM$ is either $k$, $3k/2$ or $3k$.
Here we present a few other examples (see \cite[Section 4]{HubardMochanMontero_2023_VoltageOperationsManiplexes} for details on how to see these operations as voltage operations).
\begin{eg}
    If $\cP$ is a $k$-orbit $n$-polytope, the number of orbits of the pyramid over $\cP$ is $(n+2)k/t$ where $1\leq t\leq n+2$. 
\end{eg}

\begin{eg}
    If $\cP$ is a $k$-orbit $n$-polytope, the group $\aut(\cX)\times C_2$ acts by automorphisms on the prism over $\cX$ (more details in Theorem\nobreakspace \ref {thm:LiftFixVolt}). Item\nobreakspace \ref {item:indexXY} of Proposition\nobreakspace \ref {prop:cota_orbitas} says that the number of orbits of this prism is 
    $\frac{2(n+1)k}{2t} = \frac{(n+1)k}{t}$ for some $1\leq t \leq \frac{2(n+1)}{2}=n+1$.
\end{eg}

\begin{eg}
    Let $\cX$ be a 3-maniplex (map), $\oo$  a Wythoffian operation and $A$  the set of reflections that fix a base vertex for the operation $\oo$. Then the number of orbits of $\oo(\cX)$ is given in Table\nobreakspace \ref {tab:orbits_whythoff}.

    \begin{table}
        \begin{tabular}{|c|c|c|c|c|}
        \hline
        $\oo$ & $A$ & Diagram &No. of orbits &  $t\leq$\\ \hline 
        Medial & $\{\rho_1\}$ &$\coxeter{A}{3}[2]$& $\frac{2k}{t}$ & 2\\
        Truncation & $\{\rho_0,\rho_1\}$ &$\coxeter{A}{3}[1,2]$& $\frac{3k}{t}$ & 3\\
        Truncation of the dual & $\{\rho_1,\rho_2\}$ &$\coxeter{A}{3}[2,3]$& $\frac{3k}{t}$ & 3\\
        Rhombitruncation & $\{\rho_0,\rho_2\}$ &$\coxeter{A}{3}[1,3]$& $\frac{4k}{t}$ & 4\\
         Omnitruncation & $\{\rho_0,\rho_1,\rho_2\}$ &$\coxeter{A}{3}[1,2,3]$& $\frac{6k}{t}$ & 6\\\hline
        \end{tabular}
    \caption{Number of orbits of the Wythoffian operations over a $k$-orbit map. The set $A$ is described in \cite[Section 4]{HubardMochanMontero_2023_VoltageOperationsManiplexes} and corresponds to the ringed nodes in the Coxeter diagram of the corresponding operation (in the sense of \cite[Section 5.7]{Coxeter_1973_RegularPolytopes}).}
    \label{tab:orbits_whythoff} 
    \end{table}
\end{eg}

We need to emphasise the importance of the hypothesis of $(\cY,\eta)$ preserving connectivity.
Without this hypothesis  Item\nobreakspace \ref {item:index} of Proposition\nobreakspace \ref {prop:cota_orbitas} might not hold and unusual examples arise.
If $\cY$ is any premaniplex and $\eta(W)=1$ for every $W\in\fg(\cY)$, then $\XY$ consists of $|\cX|$ disconnected copies of $\cY$. 
Its automorphism group is $\aut(\cY)\wr S_{|\cX|}$. 
The group $\aut(\cX)$ (as a subgroup of $S_{|\cX|}$) acts on $\cX \ertimes \cY$ by permuting the copies of $\cY$.
In general, the index of $\aut(\cX)$ in $\GXY$ is larger than $|\cY|$.

Moreover, even if $(\cY,\eta)$ satisfies that $\XY$ is connected for a certain premaniplex $\cX$, if the operator $(\cY,\eta)$ does not preserve connectivity, then Item\nobreakspace \ref {item:index} of Proposition\nobreakspace \ref {prop:cota_orbitas} might not hold: 

\begin{eg}[Unstable maps, see \cite{Jones_2019_UnstableMaps}]
A  non-orientable map $\cM$ is \emph{stable} if every automorphism of its orientable double cover is a lift of an automorphism of $\cM$; it is unstable otherwise. 
In \cite[Example 3.5]{HubardMochanMontero_2023_VoltageOperationsManiplexes} we showed that we can see the orientable double cover as a voltage operator $(\2_{\emptyset}, \mu)$ with  $\2_{\emptyset}$ the two-flag premaniplex with no semiedges and $\mu$ the voltage assignment that satisfies $\mu(\pth{y}) = r_{i}$ for every $0 \leq i \leq 2$ and $y \in \2_{\emptyset}$.
In our language, a stable map $\cM$ is a map that satisfies $[\aut(\cM \ertimes[\mu] \2_{\emptyset}) : \aut(\cM)] \leq 2 $.
In \cite[Example 2]{Jones_2019_UnstableMaps}, Jones shows a map $\cM$ on the Klein bottle with $|\aut(\cM)| = 8$ such that 
the orientable double cover of $\cM$ is the regular toroidal map $\left\{ 4,4 \right\}_{(2,2)}$;
the automorphism group of this map has $64$ elements. 
In this case both $\cM$ and $\cM \ertimes[\mu] \2_{\emptyset}$ are connected but Proposition\nobreakspace \ref {prop:cota_orbitas} does not hold because $(\2_{\emptyset}, \mu)$ does not preserve connectivity: if $\cK$ is an orientable map  $\cK \ertimes[\mu] \2_{\emptyset}$ consists of two isomorphic copies of $\cK$.
\end{eg}

Unfortunately, many natural operations on maps, polyhedra and polytopes do not preserve connectivity when seen as voltage operations. 
As seen above, the orientable double cover is one of them, but similar examples arise naturally when we consider operations that are only interesting for a certain family $\cF$ of maniplexes (non-orientable, non-vertex-bipartite, etc.) and give as a result an element of the complementary family $\bar{\cF}$.
Often when we apply the operation to a maniplex $\cX$ in $\bar{\cF}$ we obtain several isomorphic copies of $\cX$.

Another family of examples of voltage operations that do no preserve connectivity arise when we need to make a choice when applying the operation (e.g. the snub operation or ``truncate every second vertex''). 
Often when seen as voltage operations the result is a family of disconnected maniplexes corresponding to each of the possible choices.

\section{F-actions and voltage operations} \label{sec:Factions}
In the previous section we have characterised the automorphisms of $\cX\ertimes\cY$ that are induced by the automorphisms of $\cX$ as those fixing the $\cY$-coordinate of $\cX\ertimes\cY$.
In this section we find a necessary condition so that all automorphisms of $\cX\ertimes\cY$ are induced by those of $\cX$. 
Such condition is given in Theorem\nobreakspace \ref {teo:onlyfromX} in the language of $F$-actions. 

$F$-actions were introduced by Orbani\'c,
in \cite{Orbanic_2007_FActionsParallel},  as a tool to understand quotients of maps and hypermaps.
They can easily be also used to study maniplexes and polytopes, and, as we shall see, premaniplexes and thus voltage operations.
The ideas are similar to some of those used by Hartley in \cite{Hartley_1999_AllPolytopesAre,hartley1999more} for abstract polytopes.
In this section, we introduce $F$-actions, and refer the reader to \cite{Orbanic_2007_FActionsParallel} for details.

Let $F$ be a finitely generated and finitely presented group acting (on the left) transitively on a set $Z$, and let $z_0\in Z$ be a fixed point, that we shall call {\em root}. 
An {\em $F$-action} is a $4$-tuple $M=(\vartheta,G, Z, z_0)$, where $G$ is a group that acts (on the left) faithfully and transitively on $Z$, $\vartheta:F \to G$ is a group epimorphism such that for $z \in Z$ and $f \in F$,
 $f \cdot_F z = \vartheta(f) \cdot_G z$, and $\cdot_F, \cdot_G$ denote the (left) action the $F$ and $G$, respectively, on $Z$.

 A connected premaniplex $\cX$ can be seen as a $\cC$-action $(\vartheta, \mon(\cX), V(\cX), x_0)$, where $V(\cX)$ here denotes the set of flags of $\cX$, $x_0$ is a {\em base} flag, and $\vartheta:\cC \to \mon(\cX)$ the natural group epimorphism between $\cC$ and the monodromy group of $\cX$.

  Given two $F$-actions $M=(\vartheta_M,G_M, Z_M, z_M)$, and $N=(\vartheta_N,G_N, Z_N, z_N)$, an {\em $F$-action morphism} between them consists of a pair $(\varphi, \psi)$ such that $\varphi:Z_M \to Z_N$ is an onto mapping satisfying that $\varphi(z_M)=z_N$, and $\psi: G_M \to G_N$ is a group epimorphism satisfying that $\psi \circ \vartheta_M=\vartheta_N$ and also that for all $g\in G_M$ and $z\in Z_M$, we have $\varphi(g\cdot_M z)=\psi(g)\cdot_N \varphi(z)$, where $\cdot_M, \cdot_N$ represent the action of $G_M$ and $G_N$ on $Z_M$ and $Z_N$, respectively.
  If $(\varphi,\psi)$ are both bijections, then we have an $F$-action isomorphism.

Let $\cX$ be a connected premaniplex with base flag $x_0$, and $N:=\stab_\cC (x_0)$. Then $K=\core_\cC (N)$ is in fact the kernel $\{\omega\in \cC \mid \omega y=y \ \mathrm{for} \ \mathrm{all} \ y\in\cX\}$ of the action of $\cC$ on $\cX$, and thus the action of $\cC / H$ on $\cX$ is faithful, which says that $\mon(\cX)=\cC / K$.
By \cite[Proposition 3.1]{Orbanic_2007_FActionsParallel},  the $\cC$-action $(\vartheta, \mon(\cX), V(\cX), x)$ is isomorphic to the $\cC$-action $(\vartheta, \mon(\cX), \cC/N, N)$.
In other words, one may regard the premaniplex $\cX$ to have the flag set $\cC/N$, where the base flag $x_0$ corresponds to $N$.
In fact, this is equivalent to regard $\cX$ as the coset graph (or Schreier graph) $\Delta(\cC, N, \{r_0, r_1, \dots, r_{n-1}\})$ (see for example \cite{HubardToledo_2023_SparseGroupsNeed_preprint} for details in the case of maniplexes).

Similarly as in the case with maps and polytopes, \cite[Theorem 3.6]{Orbanic_2007_FActionsParallel} implies the following result.

\begin{proposition}\label{prop:autonorm}
    Let $\cX$ be a connected premaniplex with base flag $x_0$, and let $N=\stab_\cC(x_0)$. 
    Then $\alpha_\omega$ is an automorphism of $\cX$ satisfying that $x_0\alpha_\omega = \omega x_0 $ if and only if $\omega\in \norm_\cC (N)$, and further, $\aut(\cX) \cong \norm_\cC (N)/N$.
\end{proposition}

\subsection{$\cC$-actions and voltage operations}
Throughout this section, let $\cX$ be a connected $n$-premaniplex with base flag $x_0$ and $N:=\stab_{\cC^n} (x_0)$, and let $(\cY,\eta)$ be a $(n,m)$-voltage operator,  let $y_0\in\cY$ be a root of the connected $m$-premaniplex $\cY$ and $L:=\stab_{\cC^{m}}(y_0)$.
Thus, as discussed above, we regard $\cX$ and $\cY$ as the $\cC$-actions $(\vartheta_\cX, \mon(\cX), \cC^n/N, N)$ and $(\vartheta_\cY, \mon(\cY), \cC^m/L, L)$, respectively.

Of course, whenever $\cX \ertimes \cY$ is connected, it can be regarded as a $\cC$-action. 
In this case, we shall see that the stabiliser of $(x_0,y_0)$ can be written in terms of the stabiliser of $x_0$ and the morphism $\zeta: L \to \cC^n$ given by $\zeta(\omega)= \eta(\pth[\omega]{y_0})$. 

\begin{lemma}\label{lemma:stabXY}
Let $x_0$ and $y_0$ be the base flags of $\cX$ and $\cY$, respectively, so that $(x_0,y_0)$ is the base flag of $\cX \ertimes \cY$.
Then $\stab_{\cC^m} (x_0, y_0)=\zeta^{-1}(N)$.
In particular, if $\cX\ertimes\cY$ is connected it may be regarded as the $\cC^m$-action $(\vartheta_\eta, \mon(\cX\ertimes \cY), \cC^m/\zeta^{-1}(N), \zeta^{-1}(N))$.
\end{lemma}

\begin{proof}
  
  Start by noticing that for any $\omega\in L$ we have that 
  \begin{eqnarray}\label{eq:xi}
  \omega(x_0,y_0) = (\eta(\pth[\omega]{y_0}) x_0,\omega y_0) = (\zeta(\omega)x_0,\omega y_0) = (\zeta(\omega)x_0,y_0).
  \end{eqnarray}

  Now, $\omega \in \stab_{\cC^m} (x_0, y_0)$ implies that $(x_0,y_0)=\omega(x_0,y_0)= (\eta(\pth[\omega]{y_0}) x_0,\omega y_0)$, so $\omega \in L$. Hence, by (\ref{eq:xi}), $(x_0,y_0)=(\zeta(\omega)x_0, y_0)$, implying that $\zeta(\omega)\in N$ and thus $\omega \in \zeta^{-1}(N)$. 

  On the other hand, if $\omega \in \zeta^{-1}(N)$, then $\omega\in L$ (since $\zeta:L\to\cC^n$), and $\zeta(\omega)\in \zeta (\zeta^{-1}(N))\subseteq N$. Hence, $(x_0,y_0)=(\zeta(\omega)x_0,\omega y_0)$, which by (\ref{eq:xi}) equals $\omega(x_0,y_0)$, implying that $\omega \in \stab_{\cC^m}(x_0,y_0)$.

  We have proved that $\stab_{\cC^m}(x_0,y_0) = \zeta^{-1}({N})$. Again, by \cite[Proposition 3.1]{Orbanic_2007_FActionsParallel}, the connected component of $\cX \ertimes \cY$ containing the flag $(x_0,y_0)$ may be regarded as the $\cC^m$-action $(\vartheta_\eta, \mon(\cX\ertimes \cY), \cC^m/\zeta^{-1}(N), \zeta^{-1}(N))$. This concludes the proof.
\end{proof}

By combining Proposition\nobreakspace \ref {prop:autonorm} with Lemma\nobreakspace \ref {lemma:stabXY}, we have the following result.

\begin{corollary}\label{AutoOperationNormalizer}
If $\cX \ertimes \cY$ is connected, then $\alpha_\omega$ is an automorphism of $\cX \ertimes \cY$ satisfying that $(x_0,y_0))\alpha_\omega = \omega (x_0,y_0)) $ if and only if $\omega\in \norm_\cC \zeta^{-1}(N)$, and further,
   $$\aut(\cX \ertimes \cY) \cong \norm_{\cC^m} (\zeta^{-1}(N))/\zeta^{-1}(N).$$
\end{corollary}

Corollary\nobreakspace \ref {AutoOperationNormalizer} implies that all symmetries of $\cX \ertimes \cY$ come from elements of $\norm_{\cC^m} (\zeta^{-1}(N))$. 
We may use this to characterise, in terms of $\cC$-actions, when all the automorphisms of $\cX \ertimes \cY$ come from automorphisms of $\cX$.

\begin{proposition}\label{XautoFaction}
    If $\cX\ertimes \cY$ is connected and every automorphism of $\aut(\cX\ertimes \cY)$ is induced by an automorphism of $\aut(\cX)$, then $\norm_\cC(\zeta^{-1}(N))\subseteq L$.
    Moreover, if $(\cY,\eta)$ preserves connectivity, and $\norm_\cC(\zeta^{-1}(N))\subseteq L$, then every automorphism of $\aut(\cX\ertimes \cY)$ is induced by an automorphism of $\aut(\cX)$.
\end{proposition}

\begin{proof}
   For the first part, by Corollary\nobreakspace \ref {AutoOperationNormalizer}, if $\omega\in \norm_\cC(\zeta^{-1}(N))$, there exists ${\alpha}_\omega\in \aut(\cX\ertimes \cY)$ such that $\omega(x_0,y_0)=(x_0,y_0){\alpha}_\omega$.
    Since ${\alpha}_{\omega}$ is an automorphism of $\cX\ertimes \cY$ induced by an automorphism of $\cX$, we have that  $(x_0,y_0){\alpha}_\omega=(x_0\alpha,y_0)$, for some $\alpha\in \aut(\cX)$.
    Hence, $(\eta(\pth[\omega]y_0)x_0,\omega y_0)=\omega(x_0,y_0)=(x_0,y_0){\alpha}_\omega=(x_0\alpha,y_0)$, implying that $\omega y_0=y_0$, so $\omega\in L$.

    For the second part, let $\alpha \in \GXY$. 
    Since $\cX\ertimes \cY$ is connected, there exists $\omega\in\cC^m$ such that $(x_0,y_0)\alpha=\omega(x_0,y_0)$. By Corollary\nobreakspace \ref {AutoOperationNormalizer} $\omega \in \norm_{\cC^m}(\zeta^{-1}(N))\subseteq L$, and thus $\omega y_0 = y_0$.
    Now $(x_0,y_0)\alpha=\omega(x_0,y_0)=(\eta(\pth[\omega]y_0)x_0,\omega y_0) =(\eta(\pth[\omega]y_0)x_0, y_0) $, and since $(\cY,\eta)$ preserves connectivity, by Proposition\nobreakspace \ref {lemma:autX-1} every automorphism of $\aut(\cX\ertimes \cY)$ is induced by an automorphism of $\aut(\cX)$.
\end{proof}

\subsection{Avoiding extra automorphisms in connectivity preserving operatrions} \label{sec:ExtraSym}
We know that $\aut(\cX)$ acts by automorphisms of $\GXY$. 
 Moreso, in Section\nobreakspace \ref {sec:autoX} we showed that every automorphism of $\GXY$ that fixes the $\cY$-coordinate of one flag can be regarded as an automorphism of $\cX$ and, further, in the last section we gave a characterization, in terms of $\cC$-actions, of when every automorphism of $\GXY$ can be regarded as one from $\cX$.

In fact, we just showed that if we assume that $(\cY,\eta)$ preserves connectivity, then $\aut(\cX) = \GXY$ if and only if $\norm_\cC(\zeta^{-1}(N)) \subset L$. 
This implies that any symmetry of $\cX \ertimes \cY$ that does not come from $\aut(\cX)$ (if any) is induced by  an element of $\norm_\cC (\zeta^{-1}(N)) \setminus L$. 

Since we have characterised this extra symmetry for operations that preserve connectivity, it shall prove useful to translate Proposition\nobreakspace \ref {coro:connected} to the language of $\cC$-actions. 
Since $\zeta(\omega)=\eta(\pth[\omega]{y_0})$, for $\omega \in L$, then the natural translation is as follows.

\begin{corollary}\label{lemma:F-connectivity}
      The voltage operator $(\cY,\eta)$ preserves connectivity if and only if $\cY$ is connected and $\zeta$ is onto. In particular, $\zeta(\zeta^{-1}(N))=N$ whenever $(\cY,\eta)$ preserves connectivity.
\end{corollary}

Recall that if a premaniplex $\cZ$ is regarded as the $\cC$-action $(\vartheta, \mon(\cZ), \cC/K, K)$, for some $K\leq\cC$, then we abuse notation and refer to $\cZ$ as $\cC/K$.

\begin{lemma}\label{extrasymmg}
Let $\cX$ be a connected $n$-premaniplex and let $(\cY,\eta)$ be an $(n,m)$-voltage operator that preserves connectivity.
If $\upsilon \in \norm_\cC (\zeta^{-1}(N)) \setminus L$, then $\cX$ covers $\cZ_{\upsilon}$, where $\cZ_{\upsilon}=\cC^n/\zeta(L\cap L^{\upsilon})$.
\end{lemma}

\begin{proof}
Let $K_{\upsilon}:=L\cap L^{\upsilon}$. 
    By \cite[Proposition 3.2]{Orbanic_2007_FActionsParallel} we have that $\cX$ covers $\cZ_{\upsilon}$ if and only if $N\leq \zeta(K_{\upsilon})$. 
Let $\omega \in N$ and observe that, by Corollary\nobreakspace \ref {lemma:F-connectivity}, $\zeta^{-1}(\omega) \neq \emptyset$. 
    Take $\tilde{\omega}\in\zeta^{-1}(\omega )$ and note that $\zeta^{-1}(\omega )\subseteq \zeta^{-1}(N)\subseteq L$, which implies that $\tilde{\omega} \in L$.

    Now $\upsilon \in \norm_{\cC^m} (\zeta^{-1}(N))$ 
    implies that if we conjugate $\tilde{\omega}$ by $\upsilon^{-1}$, the result is again in $\zeta^{-1}(N)$, which is a subset of $L$. That is,
$\tilde{\omega}^{{\upsilon}^{-1}}\in L$.
    Moreover, $\tilde{\omega}\in L$ implies that $\tilde{\omega}^{{\upsilon}^{-1}}\in L^{{\upsilon}^{-1}}$, and therefore $\tilde{\omega}^{{\upsilon}^{-1}}\in L \cap L^{{\upsilon}^{-1}}$.
This in turn implies that $\tilde{\omega}\in L\cap L^{\upsilon} = K_{\upsilon}$, and thus $\zeta(\tilde{\omega})\in \zeta(K_{\upsilon})$.
    But $\tilde{\omega}\in\zeta^{-1}(\omega )$ implies that $\zeta(\tilde{\omega})=\omega$, and therefore $\omega \in \zeta(K_{\upsilon})$.
    \end{proof}

The following result is an immediate consequence of the lemma.

\begin{thm}\label{teo:onlyfromX}
    If for each ${\upsilon}\in\cC^m\setminus L$ we have that $\cX$ does not cover $\cZ_{\upsilon}$, then all the symmetries of $\cX \ertimes \cY$ are induced by symmetries of $\aut(\cX)$.
\end{thm}

Following Lemma\nobreakspace \ref {extrasymmg}, for each $\upsilon\in\cC^m$ we set $K_{\upsilon}:=L\cap L^{\upsilon}$, and $\cZ_{\upsilon}:=\cC/\zeta(K_{\upsilon})$.
By definition, there is a bijection between the sets $\{K_{\upsilon} \mid {\upsilon}\in\cC\}$ and $\{\cZ_{\upsilon}\mid {\upsilon}\in\cC\}$.
Moreover, there is a surjective function from $\{L^{\upsilon} \mid {\upsilon}\in\cC^m\}$ to $\{K_{\upsilon} \mid {\upsilon}\in\cC^m\}$ that, in general, is not a bijection.
This implies that \[|\{\cZ_{\upsilon}\mid {\upsilon}\in\cC\}| \leq |\{L^{\upsilon} \mid {\upsilon}\in\cC\}| = [\cC : \norm_{\cC}(L)]. \]
Now observe that the cosets of $\norm_{\cC}(L)$ are in correspondence with the orbits of $\aut(\cY)$.
In fact, by Proposition\nobreakspace \ref {prop:autonorm}: 
\[(\upsilon \norm_{\cC}(L)) y_{0} = \upsilon( \norm_{\cC}(L) y_{0}) = \upsilon( y_{0} \aut(\cY) ) = (\upsilon y_{0}) \aut(\cY).
\]
It follows that $[\cC : \norm_{\cC}(L)]$ coincides  with the number of orbits of  $\cY$ under the action of $\aut(\cY)$.
Therefore we have the following result.

\begin{lemma}\label{lemma:fewg}
    The number of premaniplexes $\cZ_{\upsilon}$ is bounded by the number of flag orbits of $\cY$ under $\aut(\cY)$.
\end{lemma}

\begin{figure}
    \centering
    \def\svgwidth{0.6\textwidth}
    \begingroup \makeatletter \providecommand\color[2][]{\errmessage{(Inkscape) Color is used for the text in Inkscape, but the package 'color.sty' is not loaded}\renewcommand\color[2][]{}}\providecommand\transparent[1]{\errmessage{(Inkscape) Transparency is used (non-zero) for the text in Inkscape, but the package 'transparent.sty' is not loaded}\renewcommand\transparent[1]{}}\providecommand\rotatebox[2]{#2}\newcommand*\fsize{\dimexpr\f@size pt\relax}\newcommand*\lineheight[1]{\fontsize{\fsize}{#1\fsize}\selectfont}\ifx\svgwidth\undefined \setlength{\unitlength}{650.15865776bp}\ifx\svgscale\undefined \relax \else \setlength{\unitlength}{\unitlength * \real{\svgscale}}\fi \else \setlength{\unitlength}{\svgwidth}\fi \global\let\svgwidth\undefined \global\let\svgscale\undefined \makeatother \begin{picture}(1,0.32891656)\lineheight{1}\setlength\tabcolsep{0pt}\put(0,0){\includegraphics[width=\unitlength,page=1]{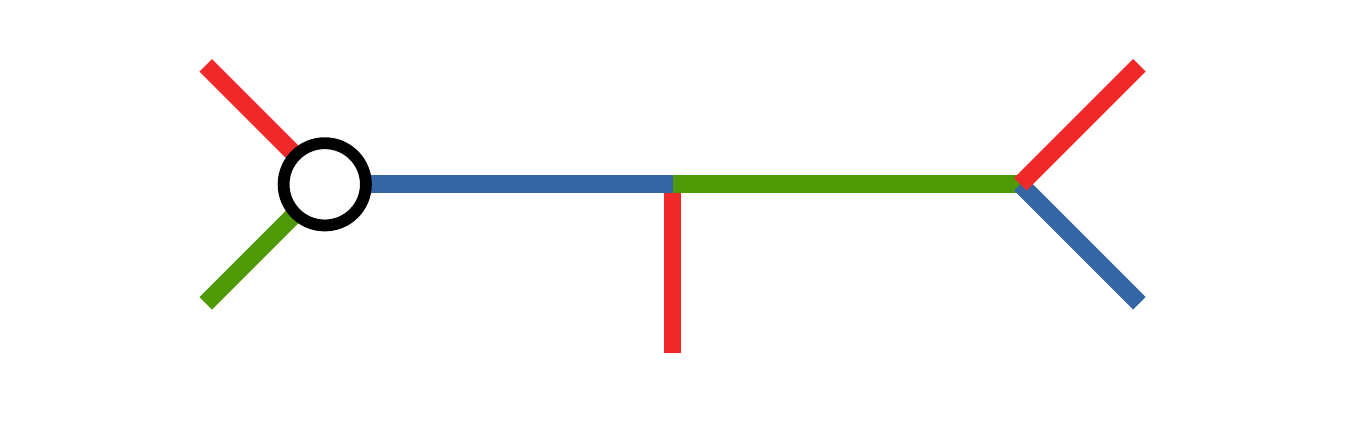}}\put(0.49585151,0.00813349){\color[rgb]{0,0,0}\makebox(0,0)[t]{\lineheight{1.25}\smash{\begin{tabular}[t]{c}$r_1$\end{tabular}}}}\put(0,0){\includegraphics[width=\unitlength,page=2]{trunc_root.pdf}}\put(0.86616887,0.29340851){\color[rgb]{0,0,0}\makebox(0,0)[lt]{\lineheight{1.25}\smash{\begin{tabular}[t]{l}$r_0$\end{tabular}}}}\put(0.13353818,0.29340851){\color[rgb]{0,0,0}\makebox(0,0)[rt]{\lineheight{1.25}\smash{\begin{tabular}[t]{r}$r_1$\end{tabular}}}}\put(0.86616887,0.06386841){\color[rgb]{0,0,0}\makebox(0,0)[lt]{\lineheight{1.25}\smash{\begin{tabular}[t]{l}$r_2$\end{tabular}}}}\put(0.13353818,0.06386841){\color[rgb]{0,0,0}\makebox(0,0)[rt]{\lineheight{1.25}\smash{\begin{tabular}[t]{r}$r_2$\end{tabular}}}}\put(0.75405463,0.18377445){\color[rgb]{0,0,0}\makebox(0,0)[t]{\lineheight{1.25}\smash{\begin{tabular}[t]{c}$y_0$\end{tabular}}}}\put(0,0){\includegraphics[width=\unitlength,page=3]{trunc_root.pdf}}\end{picture}\endgroup      \caption{The truncation operator}
    \label{fig:truncado}
\end{figure}

\begin{eg}\label{eg:trunc-y-med}
As we pointed out in the introduction, the truncation is one of the classical operations, and was studied in \cite{OrbanicPellicerWeiss_2010_MapOperations$k$}. In particular, we observe that
\cite[Proposition 4.4]{OrbanicPellicerWeiss_2010_MapOperations$k$} can be interpreted as follows: If the truncation of a map $\cM$ has more symmetries than $\cM$, then $\cM$ must cover the premaniplex $\2_{\{0,1\}}$ (the  premaniplex in Figure\nobreakspace \ref {fig:medial}).
To see this proposition as a particular case of Theorem\nobreakspace \ref {teo:onlyfromX}, consider the truncation operator ($\cY,\eta)$ described in Figure\nobreakspace \ref {fig:truncado}, together with a base flag $y_0$.
Let $\cM$ be a map (i.e. a  3-maniplex).
If the truncation of $\cM$ has more symmetries than $\cM$, Theorem\nobreakspace \ref {teo:onlyfromX} and Lemma\nobreakspace \ref {lemma:fewg} tell us that it must cover either $\cZ_{r_1}$ or $\cZ_{r_2 r_1}$. As one can see in \cite{OrbanicPellicerWeiss_2010_MapOperations$k$}, $L\cap L^{r_1} = L\cap L^{r_2 r_1} = \gen{r_0, r_0^{r_1}, r_0^{r_1 r_2}}$.
By calculating the voltages of the paths based on $y_0$ induced by these generators we see that $\zeta(L\cap L^{r_1}) = \zeta(L\cap L^{r_2 r_1}) = \gen{r_0,r_1,r_1^{r_2}}$, which is exactly the stabiliser of any vertex of $\2_{\{0,1\}}$.
Therefore, $\cZ_{r_1} = \cZ_{r_2 r_1} = \2_{\{0,1\}}$.
\end{eg}

Note that to use Theorem\nobreakspace \ref {teo:onlyfromX}, one has to check that for each ${\upsilon}\in\cC^m$, $\cX$ does not cover $\cZ_{\upsilon}$. 
This might seem like a lot of work, but Lemma\nobreakspace \ref {lemma:fewg} bounds the number of tests needed by the flag orbits of $\cY$ under $\aut(\cY)$, which is itself bounded by the number of flags of $\cY$. 
This means that if our operators have few flags, the number of tests is small.

Moreover, if there is some $\omega\in L\setminus L^\upsilon$, then $(L\cap L^\upsilon)^\omega = L \cap L^{\omega\upsilon}$, which in turn implies that $\zeta(L\cap L^\upsilon)$ and $\zeta(L \cap L^{\omega\upsilon})$ are conjugates.
Then, the main results of \cite{Hartley_1999_AllPolytopesAre} say that $\cZ_\upsilon$ and $\cZ_{\omega\upsilon}$ are isomorphic.
In terms of premaniplexes, this means that if identifying $y_0$ with $y=\upsilon y_0$ would force us to identify $y_0$ with $y_1=\omega\upsilon y_0$, then the premaniplexes $\cZ_\upsilon$ and $\cZ_{\omega\upsilon}$ are isomorphic.
Using this, we get an even smaller number of $\cZ_\upsilon$ to check, for example, if the only quotient of $\cY$ is $\1^m$, we only have to find one $\cZ_\upsilon$  ($\upsilon\notin L$) and all others will be isomorphic.

However, there is an inconvenient with Theorem\nobreakspace \ref {teo:onlyfromX}: it may happen that one of the $\cZ_\upsilon$ is $\1^n$: the $n$-premaniplex with only one point. This premaniplex is covered by any other premaniplex (of the same rank), so in this case Theorem\nobreakspace \ref {teo:onlyfromX} does not give any useful information.
This might look particularly discouraging when one notices that if $\upsilon\in\norm_{\cC^m}(L)$ (that is, when there is an automorphism of $\cY$ mapping $y_0$ to  $\upsilon y_0$), then $\cZ_\upsilon = \1^n$.
Automorphisms of $\XY$ induced by those of $\cY$ are studied in the next section.

\section{Automorphisms induced by $\cY$}\label{sec:autoY}

We have discussed in Section\nobreakspace \ref {sec:Factions} that Theorem\nobreakspace \ref {teo:onlyfromX} does not help us avoid automorphisms of $\XY$ coming from $\cY$.
These automorphisms must be studied on their own, so we dedicate this section to them. 
Later on we reformulate Theorem\nobreakspace \ref {teo:onlyfromX} considering the analysis of this section (see Theorem\nobreakspace \ref {teo:onlyfromXorY}).

Symmetries coming from $\cY$ can be divided into two classes: those that always occur, and those that may or may not occur depending on $\cX$. We will start by studying the first class.

\begin{figure}
    \centering
    \begin{scriptsize}
\def\svgwidth{\textwidth}
	\begingroup \makeatletter \providecommand\color[2][]{\errmessage{(Inkscape) Color is used for the text in Inkscape, but the package 'color.sty' is not loaded}\renewcommand\color[2][]{}}\providecommand\transparent[1]{\errmessage{(Inkscape) Transparency is used (non-zero) for the text in Inkscape, but the package 'transparent.sty' is not loaded}\renewcommand\transparent[1]{}}\providecommand\rotatebox[2]{#2}\newcommand*\fsize{\dimexpr\f@size pt\relax}\newcommand*\lineheight[1]{\fontsize{\fsize}{#1\fsize}\selectfont}\ifx\svgwidth\undefined \setlength{\unitlength}{1024.52034057bp}\ifx\svgscale\undefined \relax \else \setlength{\unitlength}{\unitlength * \real{\svgscale}}\fi \else \setlength{\unitlength}{\svgwidth}\fi \global\let\svgwidth\undefined \global\let\svgscale\undefined \makeatother \begin{picture}(1,0.52534591)\lineheight{1}\setlength\tabcolsep{0pt}\put(0,0){\includegraphics[width=\unitlength,page=1]{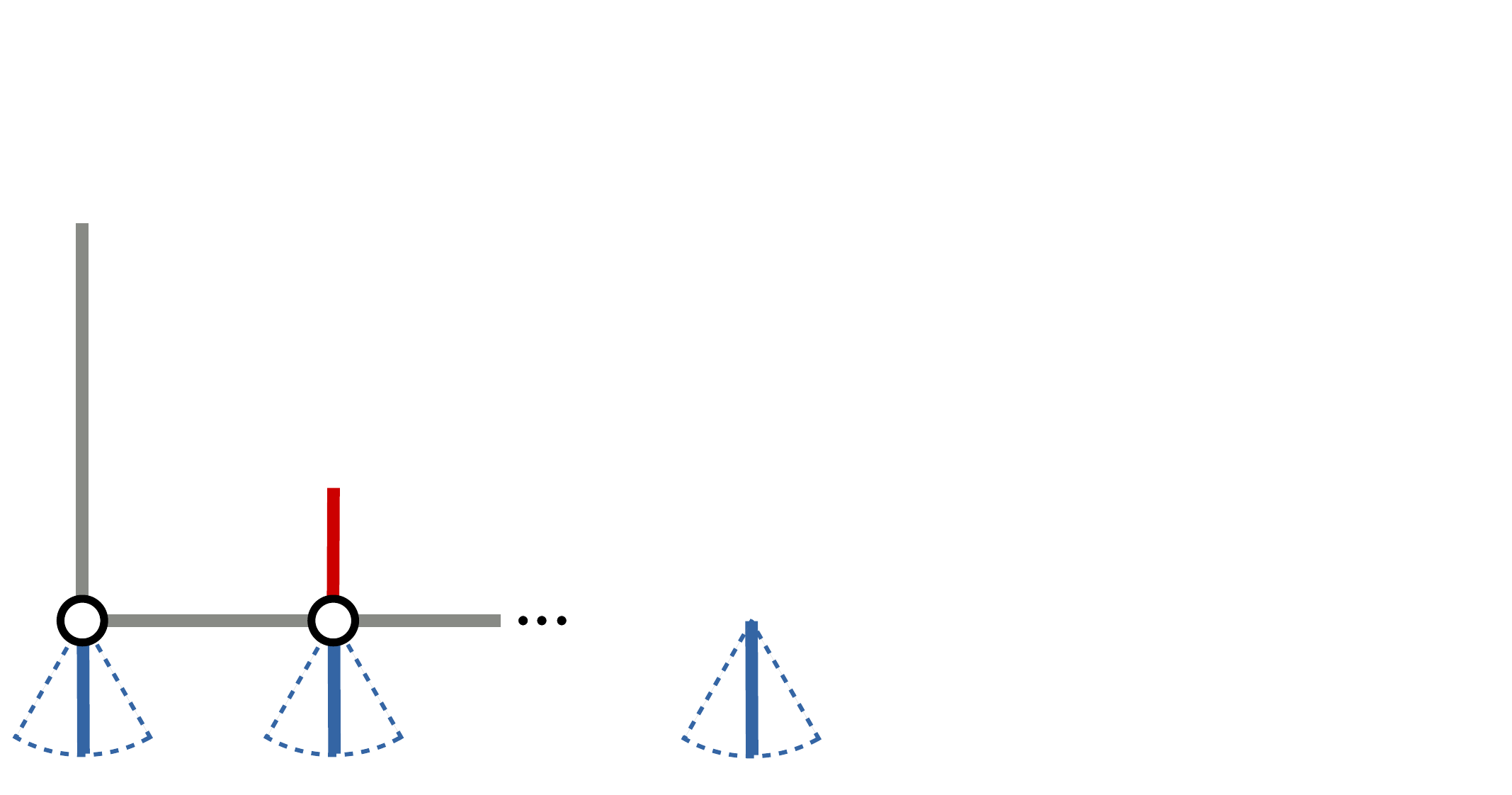}}\put(0.13748748,0.12881283){\color[rgb]{0,0,0}\makebox(0,0)[t]{\lineheight{1.25}\smash{\begin{tabular}[t]{c}$1$\end{tabular}}}}\put(0.0545253,0.00371382){\color[rgb]{0,0,0}\makebox(0,0)[t]{\lineheight{1.25}\smash{\begin{tabular}[t]{c}$r_{i-1}$\end{tabular}}}}\put(0,0){\includegraphics[width=\unitlength,page=2]{pri.pdf}}\put(0.45566975,0.12881283){\color[rgb]{0,0,0}\makebox(0,0)[rt]{\lineheight{1.25}\smash{\begin{tabular}[t]{r}$t$\end{tabular}}}}\put(0.53867384,0.12881283){\color[rgb]{0,0,0}\makebox(0,0)[lt]{\lineheight{1.25}\smash{\begin{tabular}[t]{l}$t+1$\end{tabular}}}}\put(0.22053344,0.00371373){\color[rgb]{0,0,0}\makebox(0,0)[t]{\lineheight{1.25}\smash{\begin{tabular}[t]{c}$r_{i-1}$\end{tabular}}}}\put(0.28930328,0.12881283){\color[rgb]{0,0,0}\makebox(0,0)[t]{\lineheight{1.25}\smash{\begin{tabular}[t]{c}$2$\end{tabular}}}}\put(0.49682599,0.00280235){\color[rgb]{0,0,0}\makebox(0,0)[t]{\lineheight{1.25}\smash{\begin{tabular}[t]{c}$r_{i-1}$\end{tabular}}}}\put(0.73235007,0.12881283){\color[rgb]{0,0,0}\makebox(0,0)[rt]{\lineheight{1.25}\smash{\begin{tabular}[t]{r}$n-1$\end{tabular}}}}\put(0.85685619,0.12881283){\color[rgb]{0,0,0}\makebox(0,0)[t]{\lineheight{1.25}\smash{\begin{tabular}[t]{c}$n$\end{tabular}}}}\put(0.22108158,0.22463625){\color[rgb]{0,0,0}\makebox(0,0)[t]{\lineheight{1.25}\smash{\begin{tabular}[t]{c}$r_i$\end{tabular}}}}\put(0.95376758,0.22363743){\color[rgb]{0,0,0}\makebox(0,0)[t]{\lineheight{1.25}\smash{\begin{tabular}[t]{c}$r_i$\end{tabular}}}}\put(0.4970166,0.22481723){\color[rgb]{0,0,0}\makebox(0,0)[t]{\lineheight{1.25}\smash{\begin{tabular}[t]{c}$r_i$\end{tabular}}}}\put(0.7739254,0.22363743){\color[rgb]{0,0,0}\makebox(0,0)[t]{\lineheight{1.25}\smash{\begin{tabular}[t]{c}$r_i$\end{tabular}}}}\put(0.22049171,0.253319){\color[rgb]{0,0,0}\makebox(0,0)[t]{\lineheight{1.25}\smash{\begin{tabular}[t]{c}$i < 1$\end{tabular}}}}\put(0.49642656,0.25349999){\color[rgb]{0,0,0}\makebox(0,0)[t]{\lineheight{1.25}\smash{\begin{tabular}[t]{c}$i < t $\end{tabular}}}}\put(0.77369692,0.25248526){\color[rgb]{0,0,0}\makebox(0,0)[t]{\lineheight{1.25}\smash{\begin{tabular}[t]{c}$i < n-1$\end{tabular}}}}\put(0.95353909,0.25248532){\color[rgb]{0,0,0}\makebox(0,0)[t]{\lineheight{1.25}\smash{\begin{tabular}[t]{c}$i < n$\end{tabular}}}}\put(0,0){\includegraphics[width=\unitlength,page=3]{pri.pdf}}\put(0.13746633,0.39486303){\color[rgb]{0,0,0}\makebox(0,0)[t]{\lineheight{1.25}\smash{\begin{tabular}[t]{c}$1$\end{tabular}}}}\put(0,0){\includegraphics[width=\unitlength,page=4]{pri.pdf}}\put(0.4556486,0.39486303){\color[rgb]{0,0,0}\makebox(0,0)[rt]{\lineheight{1.25}\smash{\begin{tabular}[t]{r}$t$\end{tabular}}}}\put(0.53865268,0.39486303){\color[rgb]{0,0,0}\makebox(0,0)[lt]{\lineheight{1.25}\smash{\begin{tabular}[t]{l}$t+1$\end{tabular}}}}\put(0.28928213,0.39486303){\color[rgb]{0,0,0}\makebox(0,0)[t]{\lineheight{1.25}\smash{\begin{tabular}[t]{c}$2$\end{tabular}}}}\put(0.73232887,0.39486303){\color[rgb]{0,0,0}\makebox(0,0)[rt]{\lineheight{1.25}\smash{\begin{tabular}[t]{r}$n-1$\end{tabular}}}}\put(0.856835,0.39486303){\color[rgb]{0,0,0}\makebox(0,0)[t]{\lineheight{1.25}\smash{\begin{tabular}[t]{c}$n$\end{tabular}}}}\put(0.02681553,0.23948499){\color[rgb]{0,0,0}\makebox(0,0)[rt]{\lineheight{1.25}\smash{\begin{tabular}[t]{r}$0$\end{tabular}}}}\put(0.05423344,0.51407924){\color[rgb]{0,0,0}\makebox(0,0)[t]{\lineheight{1.25}\smash{\begin{tabular}[t]{c}$r_{i-1}$\end{tabular}}}}\put(0.22024159,0.51407915){\color[rgb]{0,0,0}\makebox(0,0)[t]{\lineheight{1.25}\smash{\begin{tabular}[t]{c}$r_{i-1}$\end{tabular}}}}\put(0.49653413,0.51316777){\color[rgb]{0,0,0}\makebox(0,0)[t]{\lineheight{1.25}\smash{\begin{tabular}[t]{c}$r_{i-1}$\end{tabular}}}}\end{picture}\endgroup  	\end{scriptsize}
	\caption{Voltage operator for the prism over an $n$-polytope}\label{fig:prism}
\end{figure}

The prism over any polygon has a symmetry swapping the two lids.
More generally, the prism over any premaniplex $\cX$ admits an automorphism that swaps two facets isomorphic to $\cX$.
This symmetry can be seen as a symmetry of the prism operator itself; in Figure\nobreakspace \ref {fig:prism} this symmetry is a horizontal reflection.
Symmetries of the premaniplex $\cY$ such as the one described above often induce additional symmetry on $\cX \ertimes \cY$ regardless of the choice of $\cX$, they do exactly when they symmetry preserves the voltage assignment $\eta$.
We shall now formally explain what this means.

We say that an automorphism $\tau$ of $\cY$ \emph{lifts to $\cX\ertimes \cY$} if there is an automorphism $\tilde{\tau}$ of $\cX\ertimes \cY$ such that for every $(x,y) \in \cX \ertimes \cY$, \[(x,y)\tilde{\tau} = (x', y\tau);\] 
or equivalently, that the following diagram commutes:
\begin{equation}\label{eq:LiftsToProduct}
    \begin{tikzcd}  
        \cX\ertimes\cY \arrow[dashed]{r}{\tilde{\tau}} \arrow{d} & \cX\ertimes\cY \arrow{d}{} \\
    \cY \arrow{r}{\tau}& \cY
\end{tikzcd}
\end{equation}
where the vertical arrows are the natural projections.
The automorphism $\tilde{\tau}$ is called a \emph{lift} of $\tau$; we also say that $\tilde{\tau}$ \emph{projects} to $\tau$.
 
\begin{remark}\label{rmk:liftid}
   The automorphisms of $\cX$ acting on $\XY$ are lifts of the identity in $\cY$.
   
   In particular, if $\cX$ is connected and $(\cY,\eta)$ preserves connectivity, by Proposition\nobreakspace \ref {lemma:autX-1}, then the automorphisms of $\cX$ acting on $\XY$ coincide with the lifts of the identity in $\cY$.
\end{remark}
Of course, each $\tau\in\cY$ may have more than one lift.
In particular, if $\cX$ has non-trivial symmetries, then the identity in $\cY$ has more than one lift (actually, at least as many as elements of $\cX$).

Symmetries of $\cY$ preserving $\eta$ always lift; more precisely:

\begin{thm}\label{thm:LiftFixVolt}
 Let $(\cY,\eta)$ be a voltage operator.
 Suppose $\cY$ has an automorphism $\tau$ such that 
 \begin{equation} \label{eq:preserveVoltages}
 \eta(W\tau)=\eta(W) 
 \end{equation}
 for all paths $W\in\fg(\cY)$.
 Then, for every $n$-premaniplex $\cX$, $\tau$ has a lift  $\tilde{\tau}\in\GXY$ satisfying \[ (x,y) \tilde{\tau} = (x,y\tau)\] for every $(x,y) \in \cX \times \cY$.
\end{thm}
\begin{proof}
  Let $\tilde{\tau}$  be defined as above.
  Then, for all $i\in\{0,1,\ldots,n-1\}$,
  \[\begin{aligned}
      ((x,y)^i)\tilde{\tau} & = (\eta(\pth{y})x,y^i)\tilde{\tau}\\
      & = (\eta(\pth{y})x,(y^i)\tau)\\
      & = (\eta(\pth{y})x,(y\tau)^i)\\
      & = (\eta(\pth{y\tau})x,(y\tau)^i)\\
      & = (x,y\tau)^i\\
      & = ((x,y)\tilde{\tau})^i. 
  \end{aligned}\]
  This proves that $\tilde{\tau}$ is an automorphism of $\cX\ertimes \cY$.

\end{proof}

The proof of the above theorem implies the following result.

\begin{coro}\label{coro:productXYeta}
Let $(\cY,\eta)$ be a voltage operator and let $\aut(\cY,\eta)$ denote the subgroup of $\aut(\cY)$ consisting of all automorphisms $\tau$ of $\cY$ satisfying Equation\nobreakspace \textup {(\ref {eq:preserveVoltages})}, $ \langle \aut(\cX), \aut(\cY, \eta) \rangle =  \aut(\cX) \times \aut(\cY,\eta)$.
\end{coro}

If $\sigma\in\aut(\cX)$ and $\tau\in\aut(\cY,\eta)$, if we consider the action of $\aut(\cY,\eta)$ given in Theorem\nobreakspace \ref {thm:LiftFixVolt}, then by Corollary\nobreakspace \ref {coro:productXYeta}, $(x,y)(\sigma,\tau)=(x\sigma,y\tau)$, for $(x,y)\in\XY$.
Thus, $(\sigma,\tau)$ can be thought of as a lift of $\tau$. 
That is, every automorphism of $\cX$ gives a lift of $\tau$ (whenever $\tau\in\aut(\cY)$ lifts).

\begin{figure}
    \centering
    \def\svgwidth{0.6\textwidth}
    \begingroup \makeatletter \providecommand\color[2][]{\errmessage{(Inkscape) Color is used for the text in Inkscape, but the package 'color.sty' is not loaded}\renewcommand\color[2][]{}}\providecommand\transparent[1]{\errmessage{(Inkscape) Transparency is used (non-zero) for the text in Inkscape, but the package 'transparent.sty' is not loaded}\renewcommand\transparent[1]{}}\providecommand\rotatebox[2]{#2}\newcommand*\fsize{\dimexpr\f@size pt\relax}\newcommand*\lineheight[1]{\fontsize{\fsize}{#1\fsize}\selectfont}\ifx\svgwidth\undefined \setlength{\unitlength}{522.84320501bp}\ifx\svgscale\undefined \relax \else \setlength{\unitlength}{\unitlength * \real{\svgscale}}\fi \else \setlength{\unitlength}{\svgwidth}\fi \global\let\svgwidth\undefined \global\let\svgscale\undefined \makeatother \begin{picture}(1,0.54460132)\lineheight{1}\setlength\tabcolsep{0pt}\put(0,0){\includegraphics[width=\unitlength,page=1]{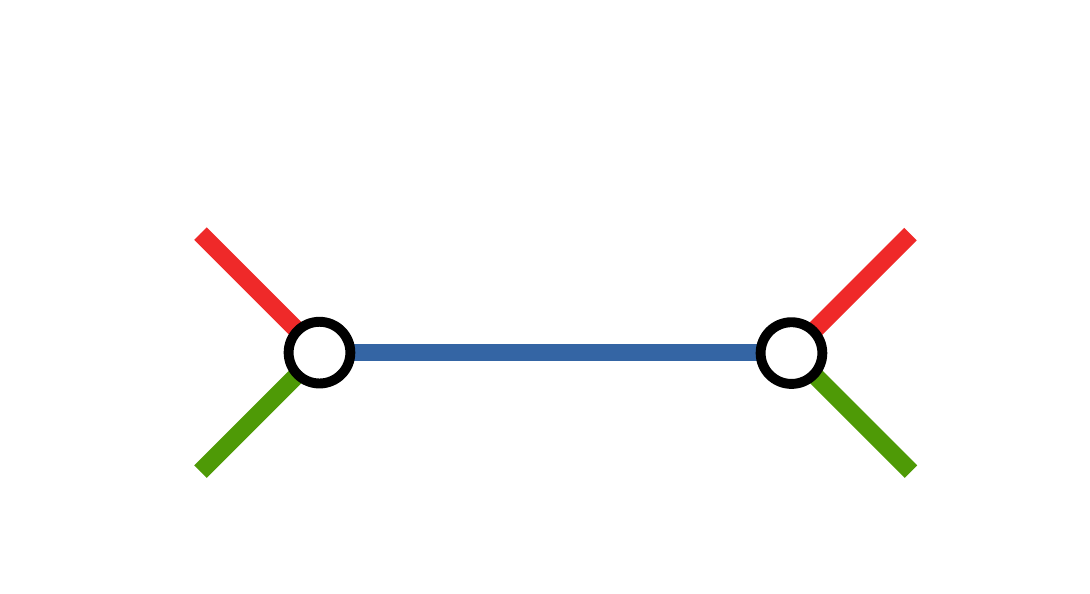}}\put(0.83049976,0.33800644){\color[rgb]{0,0,0}\makebox(0,0)[lt]{\lineheight{1.25}\smash{\begin{tabular}[t]{l}$r_1$\end{tabular}}}}\put(0.18276638,0.34319309){\color[rgb]{0,0,0}\makebox(0,0)[rt]{\lineheight{1.25}\smash{\begin{tabular}[t]{r}$r_1$\end{tabular}}}}\put(0.16605552,0.10004986){\color[rgb]{0,0,0}\makebox(0,0)[rt]{\lineheight{1.25}\smash{\begin{tabular}[t]{r}$r_0$\end{tabular}}}}\put(0.83358018,0.10499536){\color[rgb]{0,0,0}\makebox(0,0)[lt]{\lineheight{1.25}\smash{\begin{tabular}[t]{l}$r_2$\end{tabular}}}}\put(0,0){\includegraphics[width=\unitlength,page=2]{medial_sym.pdf}}\put(0.51011952,0.53135491){\color[rgb]{0,0,0}\makebox(0,0)[t]{\lineheight{1.25}\smash{\begin{tabular}[t]{c}$\tau$\end{tabular}}}}\put(0,0){\includegraphics[width=\unitlength,page=3]{medial_sym.pdf}}\end{picture}\endgroup      \caption{The medial operator}
    \label{fig:medial}
\end{figure}

Note that there might be automorphisms of $\cX\ertimes\cY$ that are not in $\aut(\cX) \times \aut(\cY,\eta)$.
For example, if $(\cY, \eta)$ is the medial operator (see Figure\nobreakspace \ref {fig:medial}), $\cY$ does not have an automorphism satisfying the conditions of the above theorem.
However, it is well-known that whenever $\cX$ is self-dual, $\cX\ertimes\cY$ has symmetries not coming from the automorphisms of $\cX$ (see \cite[Theorem 4.1]{OrbanicPellicerWeiss_2010_MapOperations$k$}).
Such extra symmetry can be though as an automorphism $\tau$ of $\cY$ that does not preserve the voltages.
Such symmetries do not necessarily lift for every premaniplex, 
Now we want to explore these cases, when the automorphisms $\tau$ of the operators lift only for some premaniplexes $\cX$, just as in the example of the medial operator.

Let $(\cY,\eta)$ be an $(n,m)$-voltage operator. Suppose that there exists $\tau\in\aut(\cY)$ such that there exists a group homomorphism $\tau^\#: \cC^n \to \cC^n$ satisfying that 
\begin{equation} \label{eq:gato}
\tau^\#(\eta(W)) = \eta(W\tau),
\end{equation}
for every $W\in \fg(\cY)$.
We shall see in Proposition\nobreakspace \ref {thm:ExtraSymmetry} that, under certain conditions, such $\tau$ lifts to an automorphism of $\cX\ertimes\cY$. 
We need to define some things in order to be able to understand the conditions of Proposition\nobreakspace \ref {thm:ExtraSymmetry}.

We can use $\tau^\#$ to define a new voltage operator $(\1^n,\tau^\#)$, where the semi-edge of colour $i$ of $\1^n$ has voltage $\tau^\#(r_i)$ for all $i\in\{0,1,\ldots,n-1\}$.
Note that by doing this we are abusing notation: $\tau^\#$ is a homomorphism, and we are also using $\tau^\#$ for the voltage assignment.
Observe that whenever $\tau^\#$ is an automorphism of $\cC^n$ this new voltage operation can be though as a $d$-automorphism (see \cite[Example 3.1]{HubardMochanMontero_2023_VoltageOperationsManiplexes} and \cite{James_1988_ComplexesCoxeterGroups—Operations}).

Now, given a premaniplex $\cX$ we define the premaniplex $\cX^{\tau^\#}$ as the premaniplex $\cX\ertimes[\tau^\#] \1^n$.

\begin{thm}\label{thm:SameResult}
    Let $(\cY,\eta)$ be an $(n,m)$-voltage operator and let $\tau\in\aut(\cY)$ be such that $\tau^\# : \cC^{n} \to \cC^{n}$ satisfies Equation\nobreakspace \textup {(\ref {eq:gato})} exists.
    Then, for every $n$-premaniplex $\cX$, the premaniplex $\cX^{\tau^\#}\ertimes \cY$ is isomorphic to $\cX \ertimes \cY$.
\end{thm}
\begin{proof}
  Using \cite[Theorem 6.1]{HubardMochanMontero_2023_VoltageOperationsManiplexes} we get that:
  $$
      \cX^{\tau^\#}\ertimes \cY = (\cX\ertimes[\tau^\#]\1^n)\ertimes \cY
      =\cX \ertimes[\theta] (\1^n \ertimes \cY),
$$
where, if $z_0$ is the only vertex of $\1^n$, 
$\theta(\pth[\omega]{(z_0,y)})=\tau^\# (\pth[{\eta(\pth[w]{y})}]{z_0})$.
Note that by definition, $\1^n \ertimes \cY \cong \cY$. This means that one can think that $\theta$ is actually a voltage assignment on $\cY$ such that $\theta(\pth[\omega]{y}):=\tau^\#(\eta(\pth[\omega]{y})) = \eta(\pth[\omega]{y\tau})$. 
Therefore, $\cX^{\tau^\#}\ertimes \cY = \cX \ertimes[\theta]\cY$,
so we actually need to prove that $\cX\ertimes[\theta] \cY$ is isomorphic to $\cX\ertimes \cY$.

  We claim that the function $(x,y)\mapsto (x,y\tau)$ is an isomorphism from $\cX\ertimes[\theta]\cY$ to $\cX\ertimes \cY$.
  Since $\tau$ is an automorphism of $\cY$, it is clearly a bijection, so to prove that it is an isomorphism, we only need to prove that the mapping preserves $i$-adjacencies for all $i\in\{0,1,\ldots,n-1\}$. When we calculate the $i$-adjacent vertex of $(x,y)$ in $\cX\ertimes[\theta]\cY$ we get
  \[
    (x,y)^i = (\theta(\pth{y})x,y^i) = (\tau^\#(\eta(\pth{y}))x,y^i) = (\eta(\pth{y\tau})x,y^i),
  \]
  and this maps to
  \[
    (\eta(\pth{y\tau})x,(y^i)\tau) = (\eta(\pth{y\tau})x,(y\tau)^i)
  \]
  which is the $i$-adjacent vertex to $(x,y\tau)$ in $\cX\ertimes \cY$.
\end{proof}

 As expected, a clear example of Theorem\nobreakspace \ref {thm:SameResult} is the fact that a polyhedron and its dual have the same medial.
 We can see that the medial operator has a unique non-trivial symmetry $\tau$ (Figure\nobreakspace \ref {fig:medial}).
 This symmetry interchanges the semi-edge with voltage $r_0$ with the one with voltage $r_2$, and also interchanges the semi-edges with voltage $r_1$.
 This means that $(\1^3,\tau^\#)$ is the duality operator $(\1^3,[r_2,r_1,r_0])$.
 Theorem\nobreakspace \ref {thm:SameResult} tells us that every 3-premaniplex $\cX$ has the same medial as its dual $\cX^{\tau^\#}=\cX^*$.
 Similarly, considering \cite[Figure 4 and Figure 7, respectively]{HubardMochanMontero_2023_VoltageOperationsManiplexes}), one can see that the snub of a polyhedron and its dual are isomorphic, and the trapezotope over a polytope and its dual are isomorphic (the latter can be seen in the figure when considering the central inversion on the cube $\bZ_2^n$, which maps every vector to its logic conjugate).

 In the case when $\cX^{\tau^\#}$ is isomorphic to $\cX$ Theorem\nobreakspace \ref {thm:SameResult} does not give any new information.
 However, in this case, we see that the isomorphism between $\cX^{\tau^\#}$ and $\cX$ together with $\tau$ induces a new automorphism of the product $\cX\ertimes\cY$, as we see in the following theorem.

 \begin{proposition}\label{thm:ExtraSymmetry}
     Let $\cX$ be a premaniplex and let $(\cY,\eta)$ be a voltage operator.
     Let $\tau\in\aut(\cY)$ be such that $\tau^\#: \cC^{n} \to \cC^{n}$ satisfying Equation\nobreakspace \textup {(\ref {eq:gato})} exists, and suppose that $\cX^{\tau^\#}$ is isomorphic to $\cX$.
     Then $\tau$ lifts to $\cX\ertimes \cY$.
 \end{proposition}
 \begin{proof}
   Let $\varphi:\cX\to \cX^{\tau^\#}$ be an isomorphism. 
   Note that, due to how the monodromy of $\cX^{\tau^\#}$ is defined, this means that $(\eta(W)x) \varphi= \eta(W\tau)(x\varphi)$ for every path $W\in\fg(\cY)$ and every vertex $x$ in $\cX$.

   We claim that the mapping $(x,y)\mapsto(x\varphi,y\tau)$ is an automorphism of $\cX\ertimes\cY$.
   In fact, for all $i\in\{0,1,\ldots,n-1\}$, the vertex $(x,y)^i = (\eta(\pth{y})x,y^i)$ is mapped to 
   \[\begin{aligned}
   ((\eta(\pth{y})x)\varphi,(y^i)\tau) &= ((\eta(\pth{y})x)\varphi,(y\tau)^i)\\&=
   ((\eta(\pth{y})x)\varphi,(y\tau)^i)\\
   &= (\eta(\pth{(y\tau)})(x\varphi),(y\tau)^i).
   \end{aligned}
   \]
  Finally we observe that $(\eta(\pth{(y\tau)})(x\varphi),(y\tau)^i)$ is exactly the $i$-adjacent vertex to $(x\varphi,y\tau)$.
 \end{proof}

Note that for any $W\in\Pi(\cY)$, $\eta(W)$ is an element of $\cC^n$, and thus one can define $\overline{\eta(W)}$ as the corresponding element of $\mon(\cX)$. 
That is, for every $x\in\cX$, $\eta(W)x=\overline{\eta(W)}x$.
Observe further that whenever $\cX\cong\cX^{\tau^\#}$, the morphism $\tau^\#:\cC^n\to\cC^n$ induces a morphism $\tau_\#:\mon(\cX)\to\mon(\cX)$ that satisfies 
\begin{equation}\label{eq:gatomon}
    \tau_\#(\overline{\eta(W)}) = \overline{\eta (W\tau)},
\end{equation}
for every $W\in\Pi(\cY)$.
However, given $\tau\in\aut(\cY)$, there may be a morphism $\tau_\#$ satisfying Equation\nobreakspace \textup {(\ref {eq:gatomon})} for some $\cX$, even if $\tau^\#:\cC^n\to\cC^n$ does not exist.
In such a case, we may define $\cX_{\tau_\#}$ as the premaniplex with the same flags as $\cX$, and the $i$-adjacencies given by $x^i=\tau_\#(r_i)x$, for every $x\in\cX$ and $i\in\{0,1,\dots, n-1\}$. The proof of the following result is exactly the same as that of Proposition\nobreakspace \ref {thm:ExtraSymmetry}.

\begin{theorem}\label{teo:idamono}
    Let $\cX$ be an $n$-premaniplex and let $(\cY,\eta)$ be an $(n,m)$-voltage operator. 
    Let $\tau\in\aut(\cY)$ be such that $\tau_\#:\mon(\cX)\to\mon(\cX)$ satisfying Equation\nobreakspace \textup {(\ref {eq:gatomon})} exists. Suppose that $\cX$ is isomorphic to $\cX_{\tau_\#}$. Then $\tau$ lifts to $\cX\ertimes\cY$.
\end{theorem}

In fact, whenever we restrict ourselves to connected premaniplexes, all the automorphisms of $\cY$ that lift to automorphisms of $\cX\ertimes\cY$ are precisely the ones covered in the above theorem, as we shall see in the following result.

\begin{theorem}\label{teo:regresomono}
    Let $\cX$ be a connected $n$-premaniplex and let $(\cY,\eta)$ be an $(n,m)$-voltage operator that preserves connectivity. 
    Let $\tau\in\aut(\cY)$ be an automorphism that lifts.
    Then there exists a morphism $\tau_\#:\mon(\cX)\to\mon(\cX)$ that satisfies Equation\nobreakspace \textup {(\ref {eq:gatomon})}, and $\cX\cong \cX_{\tau_\#}$.
\end{theorem}

\begin{proof}
    Let $\tilde{\tau}$ be a lift of $\tau$. There exists a function $\theta:\cX\ertimes\cY\to\cC$ such that $(x,y)\tilde{\tau}=(\theta(x,y)x,y\tau)$, for all $(x,y)\in\cX\ertimes\cY$.

    By using the same argument as the one used at the beginning of the proof of Lemma\nobreakspace \ref {lemma:autX}, one can see that if $y,y'\in\cY$, then $\theta(x,y)x=\theta(x,y')x$, implying that the action of $\theta$ on $x$ does not depend on the $\cY$-coordinate, and thus we may denote $\theta(x,y)$ simply as $\theta(x)$. 
    We shall abuse the notation and refer to $\theta(x)$ as an element of $\cC$ or its corresponding element in $\mon(\cX)$, depending on the situation, as we shall only use $\theta(x)$ to act on elements of $\cX$.

    Let $\omega\in\cC^n$. Then, similar to the proof of Lemma\nobreakspace \ref {lemma:autX}, on the one hand, we have that
    \begin{eqnarray}
        (\omega(x,y))\tilde{\tau}=(\eta(\pth[\omega]{y})x,\omega y)\tilde{\tau} = (\theta(\eta(\pth[\omega]{y})x)\eta(\pth[\omega]{y})x, \omega y \tau). \nonumber
    \end{eqnarray}
On the other hand,
\begin{eqnarray}
    \omega ((x,y)\tilde{\tau})=\omega (\theta(x)x,y\tau) = (\eta(\pth[\omega]{y\tau})\theta(x)x, \omega y \tau). \nonumber
\end{eqnarray}

    Since $\tilde{\tau}\in\GXY$, then $(\omega(x,y))\tilde{\tau}=\omega ((x,y)\tilde{\tau})$, implying that
    $$
        \theta(\eta(\pth[\omega]{y})x)\eta(\pth[\omega]{y})x = \eta(\pth[\omega]{y\tau})\theta(x)x,
 $$
which in turn implies that 
 \begin{equation}\label{eq:omegatau}
        \theta(\overline{\eta(\pth[\omega]{y}})x)\overline{\eta(\pth[\omega]{y}})x = \overline{\eta(\pth[\omega]{y\tau})}\theta(x)x,
    \end{equation}

    Let $\varphi:\cX \to \cX$ defined by $x\varphi=\theta(x)x$.
    We claim that $\varphi$ is a bijection. 
First note that if $x\in\cX$, then take $y\in\cY$ and observe that $\tilde{\tau}\in\GXY$, 
implies that $(\tilde{\tau})^{-1}$ is well defined. 
Hence, given $x\in\cX$, and $y\in\cY$ we have that $(x,y)(\tilde{\tau})^{-1}=(x',y')$, for some $x'\in\cX$ and $y'\in\cY$. 
This implies that $(x,y)=(x',y')\tilde{\tau}=(\theta(x')x',y'\tau)$, and so $x=\theta(x')x'$ tells us that $\varphi$ is onto.
On the other hand, if $\theta(x)x=\theta(x')x'$, then for $y\in\cY$, $(x,y)\tilde{\tau}=(\theta(x)x,y\tau)=(\theta(x')x',y\tau)=(x',y)\tilde{\tau}$. Since $\tilde{\tau}$ is a bijection, then $x=x'$.

Let us now define $\tau_\#:\mon(\cX)\to\mon(\cX)$ such that for $\omega\in\cC^n$ and $x\in\cX$,
$$
\tau_\#(\bar{\omega})x=(\bar{\omega}(x\varphi^{-1}))\varphi ;
$$
it is straightforward to see that $\tau_\#(\bar{\omega_1 \omega_2}) = \tau_\#(\bar{\omega_1})\tau_\#(\bar{\omega_2})$, so $\tau_\#$ is a morphism. 

Moreover, by definition we have that
$$
\tau_\#(\bar{\omega})(x\varphi)=(\bar{\omega}(x))\varphi = \theta(\bar{\omega}(x))\bar{\omega}x.
$$
Then, for any $\pth[\omega]{y}\in\Pi^y(\cY)$ we have that,  by Equation\nobreakspace \textup {(\ref {eq:omegatau})} ,
$$
    \tau_\#(\overline{\eta(\pth[\omega]{y})})(x\varphi)  = \theta(\overline{\eta(\pth[\omega]{y})}x)\overline{\eta(\pth[\omega]{y})}x = \overline{\eta(\pth[\omega]{y\tau})}\theta(x)x = \overline{\eta(\pth[\omega]{y\tau})} (x \varphi).
$$
Since $\varphi$ is a bijection, then $ \tau_\#(\overline{\eta(\pth[\omega]{y})})$ and $\overline{\eta(\pth[\omega]{y\tau})}$ are elements of $\mon(\cX)$ that act in the same way in every element of $\cX$, implying that 
$$
\tau_\#(\overline{\eta(\pth[\omega]{y})}) = \overline{\eta(\pth[\omega]{y\tau})},
$$
which in turn implies that $\tau_\#$ is a morphism that satisfies Equation\nobreakspace \textup {(\ref {eq:gatomon})}.
The isomorphism between $\cX$ and $\cX_{\tau_\#}$ is given by $\varphi$.

\end{proof}

A natural question to ask after the last two results is:
\begin{question}
    If $\tau\in\cY$ induces an automorphism $\tau_\#$ of $\mon(\cX)$ satisfying Equation\nobreakspace \textup {(\ref {eq:gatomon})}, for some premaniplex $\cX$, is $\cX$ isomorphic to $\cX_{\tau_\#}$?
\end{question}

Theorems\nobreakspace \ref {teo:idamono} and\nobreakspace  \ref {teo:regresomono} offer a characterisation of when an automorphism of $\cY$ lifts to an automorphism of $\cX \ertimes \cY$.
A direct consequence is that the set $\Gamma$ of automorphisms of $\cY$ that lift to $\cX \ertimes \cY$ is a subgroup of $\aut(\cY)$. 
Furthermore, if $\tilde{\Gamma}$ denotes the set of automorphisms of $\cX \ertimes \cY$ that project to $\Gamma$, then $\tilde{\Gamma}$ is a group and there is a surjective group homomorphism $\pi: \tilde{\Gamma} \to \Gamma$.  

On the other hand, if $\cX$ is connected and $(\cY, \eta)$ preserves connectivity, Proposition\nobreakspace \ref {lemma:autX-1} characterises the automorphisms of $\cX$ as the automorphisms of $\cX \ertimes \cY$ that project to $id_\cY$. 
In particular, we can think of $\Gamma(\cX) $ as a subgroup of $\tilde{\Gamma}$.
In other words, we have proved the following result.

\begin{prop}\label{prop:extension}
Let $\cX$ be a connected premaniplex, and $(\cY,\eta)$ be a voltage operator that preserves connectivity. 
Let $\Gamma \leq \Gamma(\cY)$ the group of automorphisms that lift and let $\tilde{\Gamma} \leq \GXY$ the subgroup of automorphisms of $\cX \ertimes  \cY$ that project to $\Gamma$. 
Then $\tilde{\Gamma}$ is a group extension of $\aut(\cX)$ by $\Gamma$.
That is, there exists a sequence of homomorphisms 
\[
\begin{tikzcd}
 \aut(\cX) \arrow[r, hook, "\iota"] & \tilde{\Gamma} \arrow[r, "\pi", two heads] & \Gamma 
\end{tikzcd}
\]
Such that $\iota(\aut(\cX)) = \ker(\pi) $. 

In particular, it follows that $\iota(\aut(\cX) )\triangleleft \tilde{\Gamma}$, $\tilde{\Gamma} / \iota(\aut(\cX)) \cong \Gamma$ and if $\cX$ and $\cY$ are finite, then $|\tilde{\Gamma}| = |\aut(\cX)||\Gamma|$.
\end{prop}

 As expected, the medial operator is an example of Proposition\nobreakspace \ref {thm:ExtraSymmetry}: 
 if $\tau$ is the non-trivial symmetry of $\cY$,
 as discussed above, the voltage operator $(\1^3,\tau^\#)$ is the dual operator in rank 3.
That is, $\cX^{\tau^\#}$ is the dual of $\cX$. 
Hence, the theorem says that when $\cX$ is self-dual, then the symmetry $\tau$ of $\cY$ lifts to $\tilde{\tau}\in\GXY$. 
Note that such a lift sends a flag $(x,y)$ of $\cX\ertimes\cY$ to a flag $(x\delta, y \tau)$, where $\delta$ is a duality. 
In fact, in this case, $\GXY$ has no other elements than those coming from $\aut(\cX)$ and the lifts of $\tau$. Since $\aut(\cX)$ always has index 2 in $\GXY$, it is normal.

Whenever $\cX$ is a self-dual premaniplex such that there exists a duality $\delta$ of order $2$, then $\GXY\cong\aut(\cX)\rtimes\langle\delta\rangle\cong\aut(\cX)\rtimes C_2$.
On the other hand, if $\cX$ does not have a duality of order $2$ (see~\cite{Jendrol_1989_NoninvolutorySelfDuality}), then the automorphims that do not come from $\aut(\cX)$ are lifts of $\tau$, and thus send $(x,y)$ to $(x\delta, y \tau)$, where $\delta$ is a duality. Since no duality has order two, no lift of $\tau$ has order two, implying that $\GXY$ cannot be seen as a semi-direct product of $\aut(\cX)$ by $C_2$.
In particular, this implies that the extension of Proposition\nobreakspace \ref {prop:extension} does not necessarily split. 

As pointed out before, Theorems\nobreakspace \ref {teo:idamono} and\nobreakspace  \ref {teo:regresomono} characterise the automorphisms of $\cY$ that lift to an automorphism of $\XY$. In fact, when we ask for connectivity the result is stronger, as we shall see, in the sense that
if an automorphism of $\XY$ that acts like a lift of $\tau\in\aut(\cY)$ in one flag, must actually be a lift of $\tau$.

\begin{lemma}\label{lemma:autY-1}
    Let $\XY$ be connected. Let $\gamma\in\GXY$ and $\tau\in \aut(\cY)$. If for some $(x_0,y_0)\in\XY$ and $x_1\in\cX$ it happens that $(x_0,y_0)\gamma = (x_1,y_0\tau)$, then $\gamma$ is a lift of $\tau$.
\end{lemma}
\begin{proof}
Let $(x,y)\in\XY$. Since $\XY$ is connected, there exists $\omega\in\cC$ such that $(x,y)=\omega(x_0,y_0)$. Note that in particular $\omega y_0 = y$.
    Then,
    \begin{align*}
        (x,y)\gamma &= (\omega(x_0,y_0))\gamma \\
         &= \omega ((x_0,y_0)\gamma) \\
         &= \omega (x_1,y_0\tau) \\
         &= (\eta(\pth[\omega](y_0 \tau)) x_1, \omega y_0\tau)\\
         &= (\eta(\pth[\omega](y_0 \tau)) x_1, y\tau).
    \end{align*}
\end{proof}

Note that Proposition\nobreakspace \ref {lemma:autX-1} is actually a corollary of the above lemma, in the particular case where $\tau$ is the identity automorphism of $\cY$. With this in mind, we state the following theorem, which is similar to Theorem\nobreakspace \ref {teo:onlyfromX}

\begin{theorem}\label{teo:onlyfromXorY}
    Let $\cX$ be a connected premaniplex and $(\cY,\eta)$ a voltage operator that preserves connectivity. 
    If for each ${\upsilon}\in\cC^m\setminus \norm_{\cC^m}(L)$ we find that $\cX$ does not cover $\cZ_{\upsilon}$, then every automorphism of $\XY$ is a lift of an automorphism of $\cY$. 
\end{theorem}
\begin{proof}
Let $x_0$ and $y_0$ be the base flags of $\cX$ and $\cY$, respectively.
Let $\gamma \in \GXY$, and let $\upsilon\in\cC^m$ be such that $(x_0,y_0)\gamma=\upsilon(x_0,y_0)$.
Then, by Corollary\nobreakspace \ref {AutoOperationNormalizer}, $\upsilon\in\norm_{\cC^m}(\zeta^{-1}(N))$.
    
    By hypothesis, if $\upsilon\notin \norm_{\cC^m}(L)$, then $\cX$ does not cover $\cZ_\upsilon$. 
    On the other hand, by Lemma\nobreakspace \ref {extrasymmg}, if $\upsilon\in \norm_{\cC^m}(\zeta^{-1}(N)) \setminus L$, then $\cX$ covers $\cZ_\upsilon$.
    This implies that $\upsilon \in \norm_{\cC^m}(L)$ (note that $L\subseteq \norm_{\cC^m}(L) $).

    Proposition\nobreakspace \ref {prop:autonorm} implies now that there exists $\tau\in\cY$, such that $\upsilon y_0= y_0 \tau$.
     By Lemma\nobreakspace \ref {lemma:autY-1} we have that $\gamma$ is a lift of $\tau$. (Note that in particular if  $\upsilon\in L$, then $\gamma\in\aut(\cX)$ (see Proposition\nobreakspace \ref {XautoFaction}), and thus $\gamma$ is a lift of the identity automorphism of $\cY$.)

\end{proof}

Because of Lemma\nobreakspace \ref {lemma:fewg}, the number of $\cZ_\upsilon$ to check is at most the number of orbits of $\cY$ minus one. 
 
\section{Concluding remarks}\label{sec:concluding}
    Theorem\nobreakspace \ref {teo:onlyfromXorY} tells us that, if a connected premaniplex $\cX$ does not cover the premaniplexes of a (usually small) family depending on $\cY$, all the automorphisms of $\XY$ can be described in terms of the automorphisms of $\cX$ and $\cY$.
    Using Theorem\nobreakspace \ref {teo:idamono}, Corollary\nobreakspace \ref {coro:productXYeta} and the way $\aut(\cX)$ is embedded in $\GXY$ (as described in Section\nobreakspace \ref {sec:voltageOperations}), we can summarise the description of $\sigma\in\GXY$ as a lift of an automorphism of $\cY$ as follows:

\vspace{.2cm}
\paragraph{\textbf{The general lifts}}\label{item:losDeYaveces}
The automorphism $\sigma$ is the lift of an automorphism $\tau \in \aut(\cY)$ such that there exists $\tau_{\#}:\mon(\cX) \to \mon(\cX)$ satisfying Equation\nobreakspace \textup {(\ref {eq:gatomon})},  $\varphi:\cX \to \cX_{\tau_{\#}}$ is an isomorphism and \[(x,y)\sigma = (x\varphi,y\tau)\] for every $(x,y) \in \cX \ertimes \cY$ (see Theorem\nobreakspace \ref {teo:idamono});
\vspace{.2cm}
\paragraph{\textbf{The lifts when $\tau_{\#}$  is the identity in $\mon(\cX)$}}\label{item:losDeYsiempre}

A particular case of the above is when 
 $\tau_{\#}$  is the identity in $\mon(\cX)$, then $\cX_{\tau_{\#}} = \cX$ and $\varphi$ is actually an automorphism of $\cX$. 
        In this case, $\sigma \in \aut(\cX) \times \aut(\cY, \eta)$ (see Corollary\nobreakspace \ref {coro:productXYeta});
\vspace{.2cm}
\paragraph{\textbf{The lifts induced by the identity automorphism of $\cY$}}\label{item:lasDeX}

We now see a particular case of the above situation. That is, 
if $\tau = id_{\cY}$ (which in particular implies that $\tau_{\#}$ is the identity in $\mon(\cX)$), then $\sigma \in \aut(\cX)$ (see Lemma\nobreakspace \ref {lemma:autX}).
\vspace{.2cm}

It is now evident that the additional symmetries of the medial of a self-dual polyhedron or that of the prism over a polytope that is not a prism are included in the preceding overview.
Nevertheless, for certain premaniplexes $\cX$ and operators $(\cY,\eta)$, there may be automorphisms of $\cX \ertimes \cY$ that do not fit the above description. Examples of this include the prism over hypercubes or the truncation of some maps.

 Let us go back to consider the truncation operator $(\cY,\eta)$ (see Figure\nobreakspace \ref {fig:truncado}). 
    Then $\cY$ has no non-trivial automorphisms.
    Therefore, the above discussion seems to imply that the only automorphisms of $\cX \ertimes \cY$ that we consider are those from $\cX$. 
    However, for example, the truncation of the regular tessellation of the plane by equilateral triangles is again a regular tessellation, so it actually has symmetries not coming from those of the triangular tessellation. 
     Another example of a truncation with more symmetry than the one originating from $\cX$ is that of the map on the sphere of type $\{2,4\}$:
   it is straightforward to see that if $(\cY, \eta)$ is the truncation operator, then $\cX \ertimes \cY$ is the cube $\cP = \{4,3\}$ and the symmetries of the cube induced by $\cX$ are those that do not move a pair of opposite faces (see Example\nobreakspace \ref {eg:trunc-y-med} and \cite[Propositions 4.3 and 4.4]{OrbanicPellicerWeiss_2010_MapOperations$k$}, for more details about the truncation having additional symmetry than that coming from the original map).

   Observe that the key to understanding the possible additional symmetry in the truncation is to recognise $\cX$ itself as the result of another operation: the medial operation.
   Let us explain this idea a bit more in more detail.

     \begin{figure}
    \centering
    \begin{subfigure}[b]{0.4\textwidth}
      \centering
      \includegraphics[width=\textwidth]{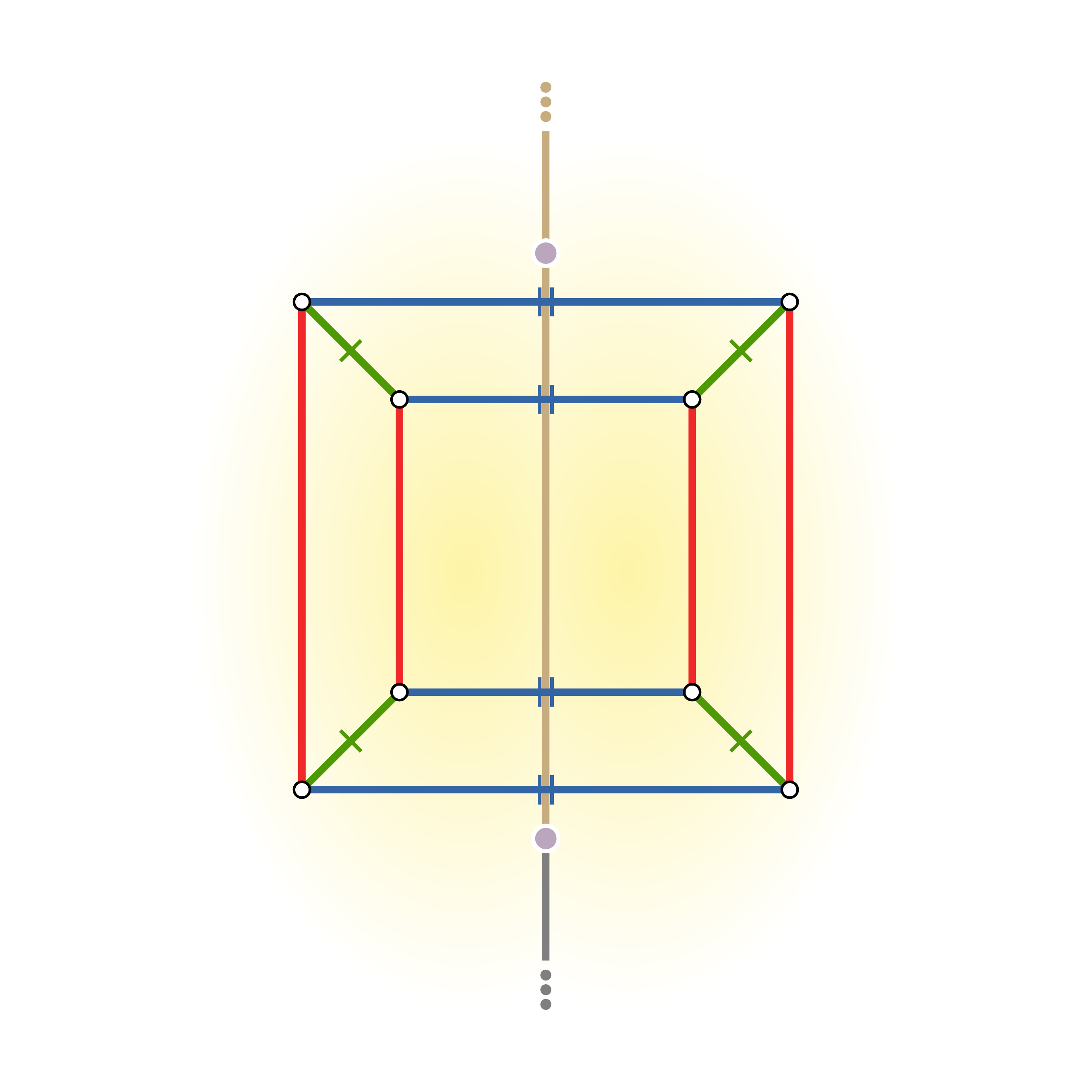}
      \caption{$\left\{ 2,2 \right\}$}
      \label{fig:22}
  \end{subfigure}
  \begin{subfigure}[b]{0.4\textwidth}
    \centering
    \includegraphics[width=\textwidth]{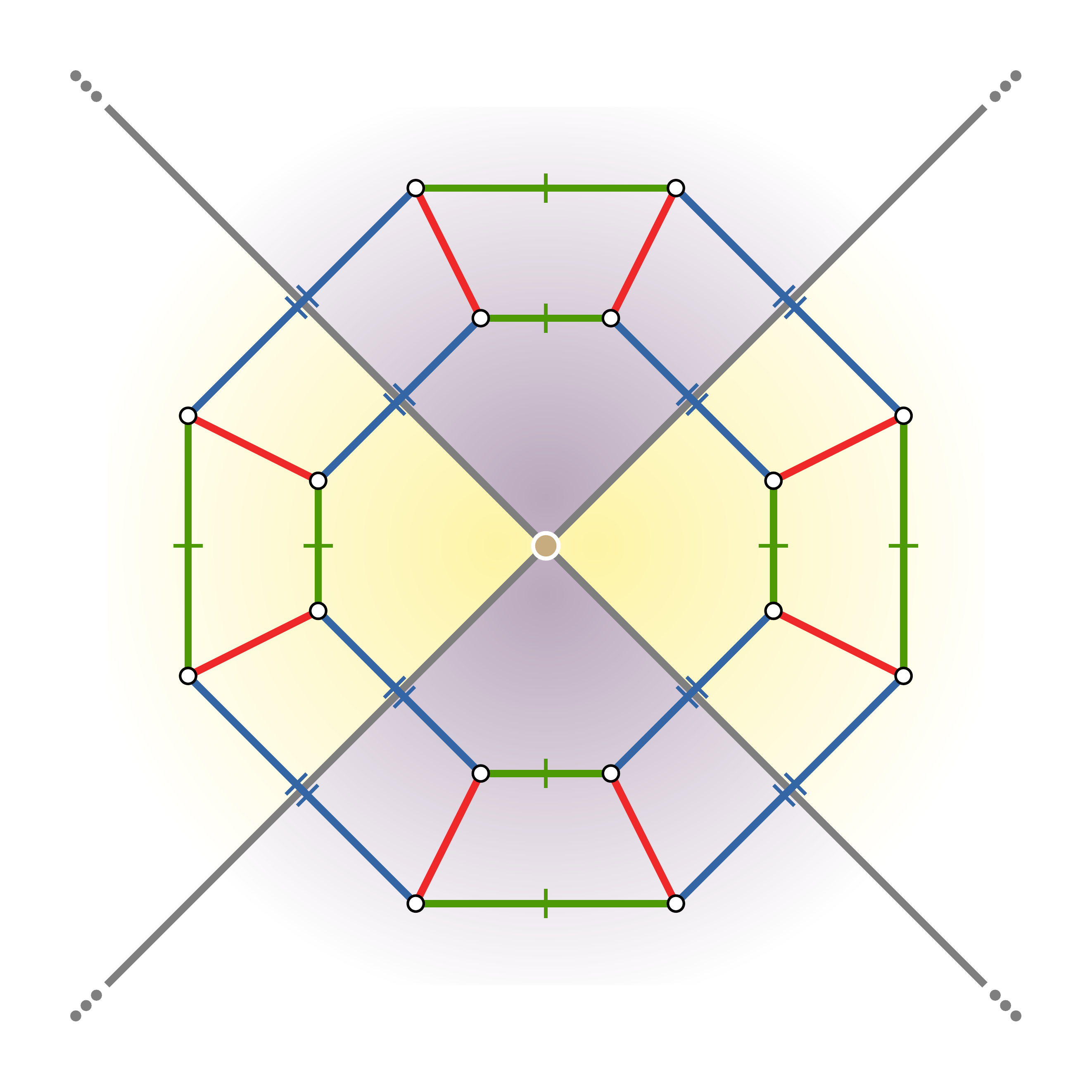}
    \caption{$\left\{ 2,4 \right\} = \med\left\{ 2,2 \right\}$}
    \label{fig:24}
\end{subfigure}
  \begin{subfigure}[b]{0.4\textwidth}
    \centering
    \includegraphics[width=\textwidth]{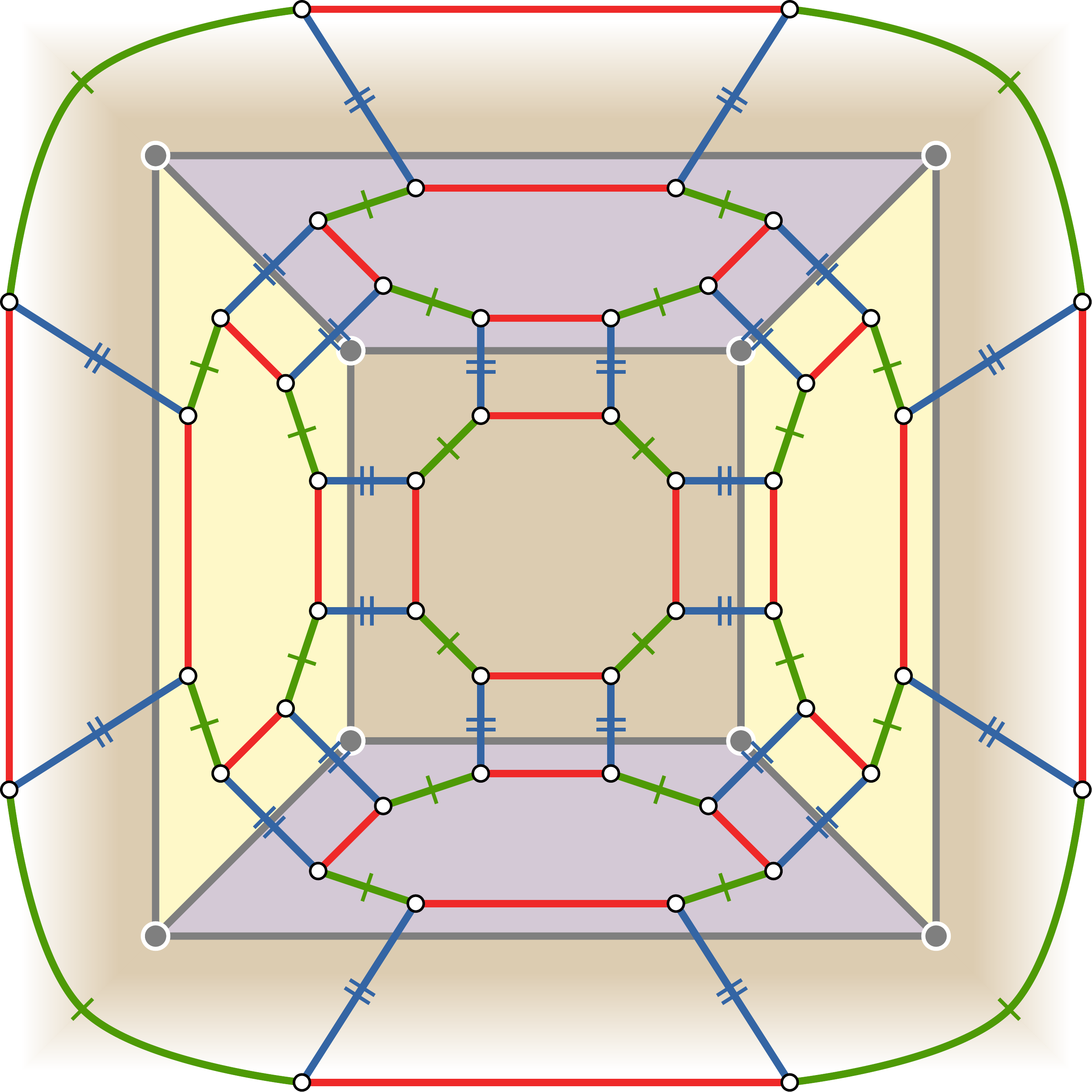}
    \caption{$\left\{ 4,3 \right\} = \tr\left\{ 2,4 \right\}$}
    \label{fig:43}
\end{subfigure}
\caption{The cube $\left\{ 4,3 \right\}$ as the truncated $\left\{ 2,4 \right\}$, which is itself the medial of $\left\{ 2,2 \right\}$.}
    \label{fig:tr(24)}
\end{figure}   

 Again, denote by $(\cY, \eta)$ the truncation operator, and let $(\wW, \vartheta)$ be the medial operator. 
 Consider $\zZ = \left\{ 2,2 \right\}$ (see Figure\nobreakspace \ref {fig:tr(24)}), and observe that $\zZ \ertimes[\vartheta] \wW$ is precisely $\cX$, the map $\{2,4\}$.
     By \cite[Theorem 6.1]{HubardMochanMontero_2023_VoltageOperationsManiplexes},  $ \cP = \cX \ertimes Y \cong  \zZ \ertimes[\theta](\wW \ertimes Y)$, where $\theta = \theta(\vartheta,\eta)$.

 \begin{figure}
    \centering
    \def\svgwidth{\textwidth}
    \begingroup \makeatletter \providecommand\color[2][]{\errmessage{(Inkscape) Color is used for the text in Inkscape, but the package 'color.sty' is not loaded}\renewcommand\color[2][]{}}\providecommand\transparent[1]{\errmessage{(Inkscape) Transparency is used (non-zero) for the text in Inkscape, but the package 'transparent.sty' is not loaded}\renewcommand\transparent[1]{}}\providecommand\rotatebox[2]{#2}\newcommand*\fsize{\dimexpr\f@size pt\relax}\newcommand*\lineheight[1]{\fontsize{\fsize}{#1\fsize}\selectfont}\ifx\svgwidth\undefined \setlength{\unitlength}{1394.170835bp}\ifx\svgscale\undefined \relax \else \setlength{\unitlength}{\unitlength * \real{\svgscale}}\fi \else \setlength{\unitlength}{\svgwidth}\fi \global\let\svgwidth\undefined \global\let\svgscale\undefined \makeatother \begin{picture}(1,0.5596891)\lineheight{1}\setlength\tabcolsep{0pt}\put(0,0){\includegraphics[width=\unitlength,page=1]{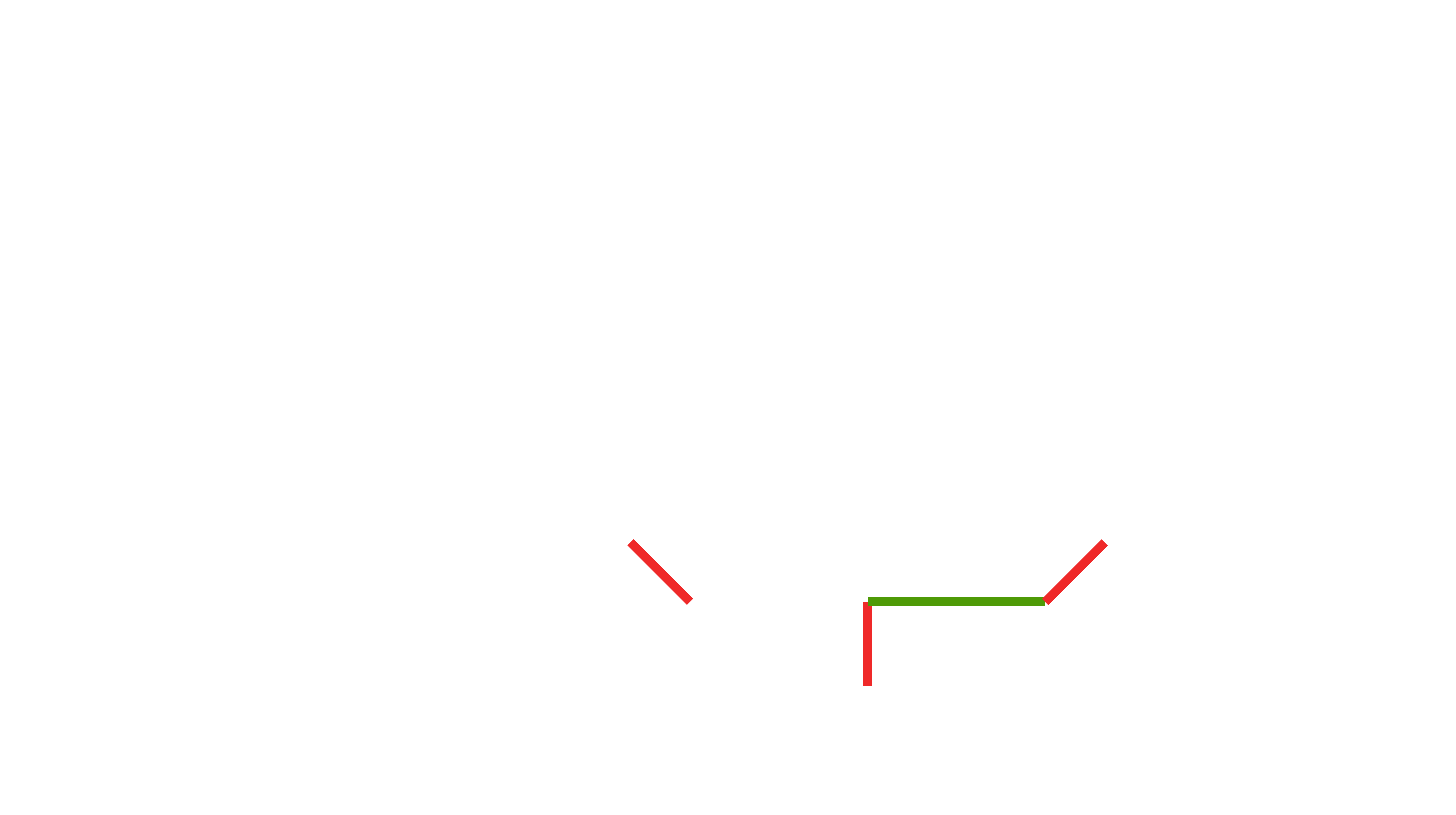}}\put(0.59597659,0.05515914){\color[rgb]{0,0,0}\makebox(0,0)[t]{\lineheight{1.25}\smash{\begin{tabular}[t]{c}$r_1$\end{tabular}}}}\put(0.7687996,0.19748405){\color[rgb]{0,0,0}\makebox(0,0)[lt]{\lineheight{1.25}\smash{\begin{tabular}[t]{l}$r_0$\end{tabular}}}}\put(0.42315346,0.19748405){\color[rgb]{0,0,0}\makebox(0,0)[rt]{\lineheight{1.25}\smash{\begin{tabular}[t]{r}$r_1$\end{tabular}}}}\put(0.7687996,0.08565733){\color[rgb]{0,0,0}\makebox(0,0)[lt]{\lineheight{1.25}\smash{\begin{tabular}[t]{l}$r_2$\end{tabular}}}}\put(0.42315346,0.08565733){\color[rgb]{0,0,0}\makebox(0,0)[rt]{\lineheight{1.25}\smash{\begin{tabular}[t]{r}$r_2$\end{tabular}}}}\put(0,0){\includegraphics[width=\unitlength,page=2]{meXtr.pdf}}\put(0.31132677,0.52279813){\color[rgb]{0,0,0}\makebox(0,0)[lt]{\lineheight{1.25}\smash{\begin{tabular}[t]{l}$r_1$\end{tabular}}}}\put(0.32149283,0.28897865){\color[rgb]{0,0,0}\makebox(0,0)[lt]{\lineheight{1.25}\smash{\begin{tabular}[t]{l}$r_1$\end{tabular}}}}\put(0.2299983,0.52279813){\color[rgb]{0,0,0}\makebox(0,0)[rt]{\lineheight{1.25}\smash{\begin{tabular}[t]{r}$r_0$\end{tabular}}}}\put(0.18933404,0.40080534){\color[rgb]{0,0,0}\makebox(0,0)[rt]{\lineheight{1.25}\smash{\begin{tabular}[t]{r}$(\cZ,\vartheta)$\end{tabular}}}}\put(0.59597659,0.00432885){\color[rgb]{0,0,0}\makebox(0,0)[t]{\lineheight{1.25}\smash{\begin{tabular}[t]{c}$(\cY,\eta)$\end{tabular}}}}\put(0.21983223,0.28897865){\color[rgb]{0,0,0}\makebox(0,0)[rt]{\lineheight{1.25}\smash{\begin{tabular}[t]{r}$r_2$\end{tabular}}}}\put(0,0){\includegraphics[width=\unitlength,page=3]{meXtr.pdf}}\put(0.42315351,0.28897865){\color[rgb]{0,0,0}\makebox(0,0)[rt]{\lineheight{1.25}\smash{\begin{tabular}[t]{r}$r_2$\end{tabular}}}}\put(0.59597658,0.26864652){\color[rgb]{0,0,0}\makebox(0,0)[t]{\lineheight{1.25}\smash{\begin{tabular}[t]{c}$r_2$\end{tabular}}}}\put(0,0){\includegraphics[width=\unitlength,page=4]{meXtr.pdf}}\put(0.42315351,0.52279813){\color[rgb]{0,0,0}\makebox(0,0)[rt]{\lineheight{1.25}\smash{\begin{tabular}[t]{r}$r_0$\end{tabular}}}}\put(0.59597658,0.54313026){\color[rgb]{0,0,0}\makebox(0,0)[t]{\lineheight{1.25}\smash{\begin{tabular}[t]{c}$r_0$\end{tabular}}}}\put(0.76879959,0.52279813){\color[rgb]{0,0,0}\makebox(0,0)[lt]{\lineheight{1.25}\smash{\begin{tabular}[t]{l}$r_1$\end{tabular}}}}\put(0.76879959,0.28897865){\color[rgb]{0,0,0}\makebox(0,0)[lt]{\lineheight{1.25}\smash{\begin{tabular}[t]{l}$r_1$\end{tabular}}}}\put(0,0){\includegraphics[width=\unitlength,page=5]{meXtr.pdf}}\end{picture}\endgroup      \caption{If  $(\wW,\vartheta)$ and $(\cY,\eta)$ are the medial and truncation operators, respectively, then  $(\wW \ertimes \cY, \theta)$ is the omnitruncation, where $\theta=(\vartheta,\eta)$.}
    \label{fig:truncadoxmedial}
\end{figure}

     Now, as seen in Figure\nobreakspace \ref {fig:truncadoxmedial}, $(\wW \ertimes \cY, \theta)$ is the omnitruncation operator \cite[Fig. 3E]{HubardMochanMontero_2023_VoltageOperationsManiplexes}. 
     The premaniplex $(\wW \ertimes \cY)$ has $S_3$ as automorphism group, and every automorphism $\tau$  induces an isomorphism $\tau_{\#}: \mon(\zZ) \to \mon(\zZ)$ (notice that $\aut(\zZ) \cong \mon(\zZ) \cong \left\langle r_{0} \right\rangle \times \left\langle r_{1} \right\rangle \times \left\langle r_{2} \right\rangle $) satisfying Equation\nobreakspace \textup {(\ref {eq:gatomon})}.
     Theorem\nobreakspace \ref {teo:idamono} implies that $\tau$ lifts to $\cP = \zZ \ertimes[\theta] (\wW \ertimes \cY) $.
     It follows from Proposition\nobreakspace \ref {prop:extension} that the group $\tilde{\Gamma} \leq (\cP)$  induced by lifts of $\aut(\wW \ertimes Y)$ is a group extension of $\aut(\zZ)$ by $S_3$. 
     In particular, $|\tilde{\Gamma}| = 48 = \aut(\cP) $, that is, every automorphism of the group $\cP$ is a lift of an automorphism of $\wW \ertimes \cY$.

\begin{figure}
    \centering
    \begin{subfigure}[b]{0.8\textwidth}
      \centering
      \includegraphics[width=\textwidth]{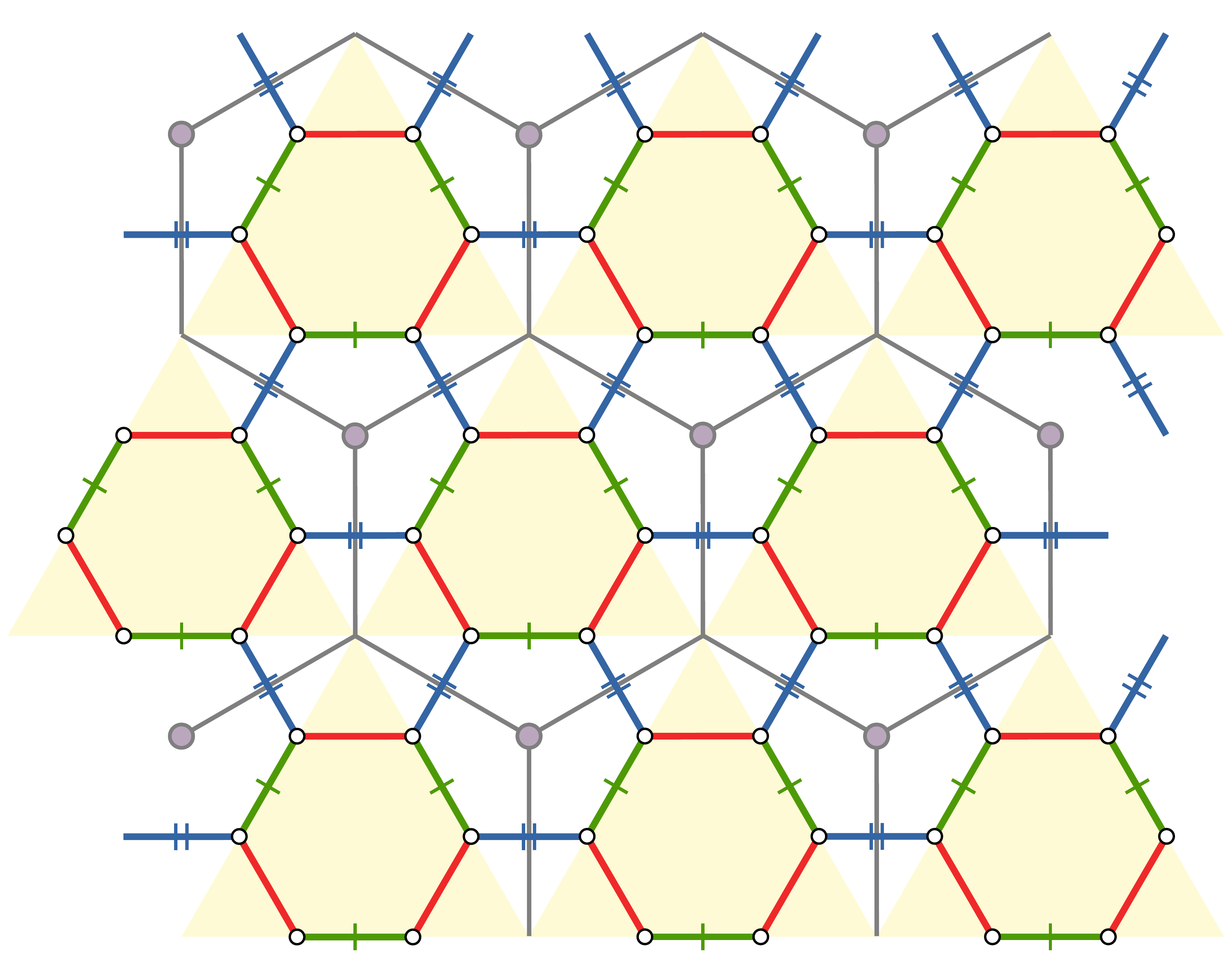}
      \caption{$(3,3,3)$}
      \label{fig:333}
  \end{subfigure}
  \begin{subfigure}[b]{0.8\textwidth}
    \centering
    \includegraphics[width=\textwidth]{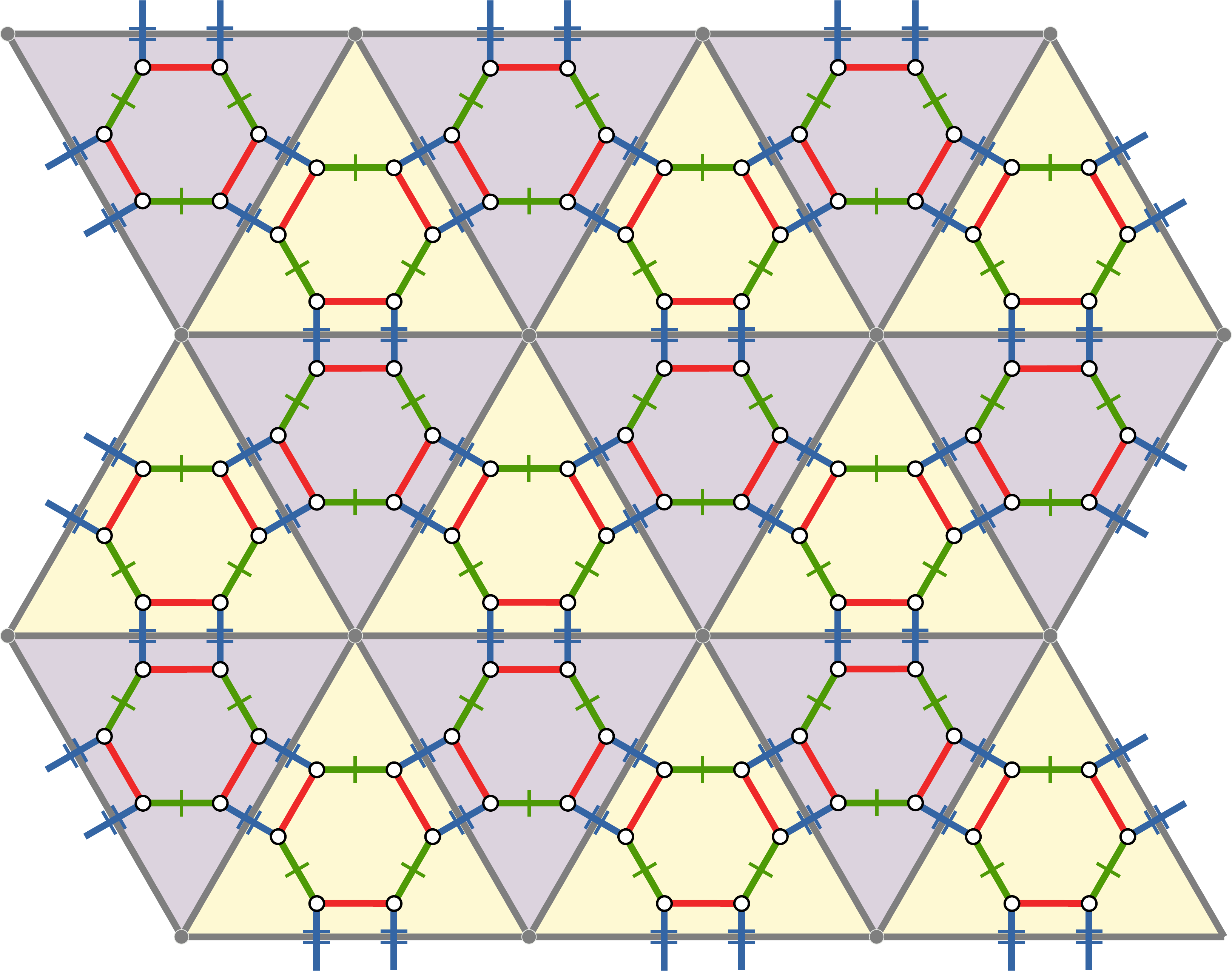}
    \caption{$\left\{ 3,6 \right\} = \med(3,3,3)$}
    \label{fig:36}
\end{subfigure}
    \caption{The tiling $\left\{ 3,6 \right\}$ can be seen as the medial of the hypermap $(3,3,3)$.}
    \label{fig:tr(36)}
\end{figure}

    One may wonder about the regular triangle tessellation and its truncation: can we understand the additional symmetry in a similar matter?
    The answer is yes!, if done with caution.
    In this case, as shown in Figure\nobreakspace \ref {fig:tr(36)}, the suitable $\zZ$ is (the flag-graph of) the regular (universal) hypermap of type $(3,3,3)$ (see \cite{CORN1988337} for details about regular hypermaps). 

In other words, if we extend our universe and drop the condition in our definition of premaniplex that if $i,j \in \left\{0, \dots, n-1 \right\} $ are such that $\left| i-j \right| \geq 2 $, then the alternating paths of length 4 with colours $i,j$ are closed, then we have {\em precomplexes} (see \cite{Wilson_2012_ManiplexesPart1} for the definition of a complex).
    In this way, we can always understand the additional symmetry of the truncation as a result of the original premaniplex being the medial of a precomplex.
    This can be seen as follows: by \cite[Proposition 4.4]{OrbanicPellicerWeiss_2010_MapOperations$k$} if the truncation of $\cX$ has extra symmetry, then $\cX$ covers the premaniplex $\2_{\{0,1\}}$ (see Example\nobreakspace \ref {eg:trunc-y-med}).
    This implies that the faces of $\cX$ can be coloured with two colours in such a way that faces that share an edge have different colours.
    Hence, we can think of the hypermap $\zZ$ whose faces correspond to the faces of one colour, its vertices to the faces of the other colour, and its hyperedges to the vertices of $\cX$.
    It is not difficult then to convince oneself that $\cX$ is the medial of $\zZ$.
    Hence, the possible additional symmetry of $\tr(\cX)$ appears if and only if $\zZ$ is self-dual.
    Note further that the premaniplex $\cZ$ actually coincides with the premaniplex $\2_{\{0,1\}}$.
    This is not coincidence, the maniplex $\2_{\{0,1\}}$ is precisely a $\cZ_\upsilon$, where $\cZ_\upsilon$ is as in Lemma\nobreakspace \ref {extrasymmg} (see Example\nobreakspace \ref {eg:trunc-y-med}).
    
   That is, in the particular case of the truncation, we have been able to understand the additional symmetry. 
   Can we do the same for any other voltage operations?
   Note that, for the prism or pyramid operators, this has been answered in \cite{GleasonHubard_2018_ProductsAbstractPolytopes}, where they are considered as specific products:
the prism (resp. pyramid) over $\cP$ has additional symmetry if and only if $\cP$ is itself a prism (resp. pyramid)(this is a particular case of \cite[Theorem A]{GleasonHubard_2018_ProductsAbstractPolytopes}). 
That is, the
additional symmetry of the prism or pyramid of a $n$-maniplex $\cX$ arises when $\cX$ can be written as the prism or pyramid of a $(n-1)$-maniplex.

As voltage operations, if $(\cY_n,\eta_n)$ denotes the prism or pyramid operator over an $n$-polytope, 
\cite[Theorem A]{GleasonHubard_2018_ProductsAbstractPolytopes} says that $\cP\ertimes[{\eta_n}] \cY_n$ has additional symmetry if and only if $\cP \cong \cW\ertimes[\eta_{n-1}] \cY_{n-1}$, for some $(n-1)$-polytope $\cW$.
In particular, $\cP$ covers $\cY_{n-1}$.
   On the other hand, Theorem\nobreakspace \ref {teo:onlyfromXorY} tells us that if $\cP\ertimes[{\eta_n}] \cY_n$ has additional symmetry, then $\cP$ covers some premaniplex $\cZ_\upsilon$ (defined as in Lemma\nobreakspace \ref {extrasymmg}). 
   For small values of $n$, with the help of GAP (\cite{gap}) one can see that, in fact, $\cY_{n-1}$ coincides with all the $\cZ_\upsilon$ coming from Theorem\nobreakspace \ref {teo:onlyfromXorY}.

   This leads us to finish this paper with the following conjecture.

   \begin{conjecture*}
       Let $\cX$ be a connected premaniplex and let $(\cY,\eta)$ be a voltage operator that preserves connectivity. 
       If $\cX$ cannot be regarded as $\cX \cong \zZ\ertimes[\vartheta] \wW_\upsilon$, for a  connectivity-preserving operator $(\wW_\upsilon,\vartheta)$, where $\cZ_\upsilon$ is defined as in Lemma\nobreakspace \ref {extrasymmg}, then all the elements of $\GXY$ are lifts of automorphisms of $\cY$.
   \end{conjecture*}
     \section*{Acknowledgements}
The authors thank the financial support of CONAHCyT grant A1-S-21678.
This paper was completed while the second author held a Post Doctoral Research Associate position in the Department of Mathematics at Northeastern University.
The third author was partially supported by ARRS Projekt 21 N1-0216.

\printbibliography
\end{document}